\documentstyle[12pt]{article}
\newlength{\abstractwidth}
\renewcommand{\thefootnote}{\fnsymbol{footnote}}
\renewcommand{\thanks}[1]{\footnote{#1}} 
\newcommand{\starttext}{ \setcounter{footnote}{0}
\renewcommand{\thefootnote}{\arabic{footnote}}}

\newcommand{\be}{\begin{equation}}
\newcommand{\bea}{\begin{eqnarray}}
\newcommand{\eea}{\end{eqnarray}} \newcommand{\ee}{\end{equation}}
 
 \def\ba{\begin{eqnarray}}
\def\ea{\end{eqnarray}}

\def\E{{\cal E}}

\def\o{\omega}
\def\Re{{\rm Re}}

\def\log{\,{\rm log}\,}

\def\o{\omega}

\def\o{\omega}

\def\na{\nabla}

\def\R{{\bf R}}

\def\p{\partial}

\def\ddb{{\partial\bar\partial}}

\def\na{{\nabla}}

\def\[{{\bf [}}
\def\]{{\bf ]}}

\begin{document}
\starttext \baselineskip=18pt \setcounter{footnote}{0}
\newtheorem{theorem}{Theorem}
\newtheorem{lemma}{Lemma}
\newtheorem{corollary}{Corollary}
\newtheorem{definition}{Definition}
\newtheorem{conjecture}{Conjecture}
\newtheorem{proposition}{Proposition}

\begin{center}
{\Large \bf FU-YAU HESSIAN EQUATIONS}

\bigskip

{\large Duong H. Phong, Sebastien Picard, and Xiangwen Zhang} \\

\medskip

\begin{abstract}

\medskip
\small{
We solve the Fu-Yau equation for arbitrary dimension and arbitrary slope 
$\alpha'$. Actually we obtain at the same time a solution of the open case 
$\alpha'>0$, an improved solution of the known case $\alpha'<0$, and solutions for a family of Hessian equations which includes the Fu-Yau equation as a special case. The method is based on the introduction of a more stringent ellipticity condition than the usual $\Gamma_k$ admissible cone condition, and which can be shown to be preserved by precise estimates with scale.
}

\end{abstract}

\end{center}

\baselineskip=15pt
\setcounter{equation}{0}
\setcounter{footnote}{0}

\section{Introduction}
\setcounter{equation}{0}

The main goal of this paper is to solve the following non-linear partial differential equation proposed in 2008 by J.X. Fu and S.T. Yau \cite{FY1},
\be \label{FY-form-eq}
i \ddb (e^u \hat{\omega} - \alpha' e^{-u} \rho) \wedge \hat{\omega}^{n-2} + \alpha' i \ddb u \wedge i \ddb u \wedge \hat{\omega}^{n-2} + \mu \, \hat{\omega}^n = 0.
\ee
Here the unknown is a scalar function $u$ on a compact $n$-dimensional K\"ahler manifold $(X,\hat\o)$, and the given data is a real $(1,1)$ form $\rho$, a function $\mu$, and a number $\alpha' \in {\bf R}$ called the slope. A key innovation in the solution is the introduction of an ellipticity condition which is more restrictive than the usual cone conditions for fully non-linear second order partial differential equations, but which can be shown to be preserved by the continuity method using some precise estimates with scale. This innovation may be useful for
other equations as well, and we shall illustrate this by using it to solve a whole family of Hessian equations in which the equation (\ref{FY-form-eq}) fits as only the simplest example.

\smallskip

The equation (\ref{FY-form-eq}) is a generalization of an equation in complex dimension $2$, which was shown in \cite{FY1} to arise from the Hull-Strominger system \cite{H1,H2, S}. The Hull-Strominger system is an extension of a proposal of Candelas, Horowitz, Strominger, and Witten \cite{CHSW} for supersymmetric compactifications of the heterotic string. It poses new geometric difficulties as it involves quadratic expressions in the curvature tensor, but it can potentially lead to a new notion of canonical metric in non-K\"ahler geometry. From our point of view, the equation (\ref{FY-form-eq}) is of particular interest as a model equation for an eventual extension of the classical theory of Monge-Amp\`ere equations of Yau \cite{Y} and Hessian equations of Caffarelli, Nirenberg, and Spruck \cite{CNS}, to more general equations mixing the unknown, its gradient, and several Hessians.

\smallskip
When the dimension of $X$ is $n=2$, the equation (\ref{FY-form-eq}) was solved by Fu and Yau in two separate papers, \cite{FY1} for the case when $\alpha'>0$, and \cite{FY2} for the case when $\alpha'<0$ (when $\alpha'=0$, the equation poses no difficulty as it reduces essentially to the Laplacian). As we shall discuss below, in the approach of \cite{FY1, FY2}, the required estimates in the two cases $\alpha'>0$ and $\alpha'<0$ are quite different. In an earlier paper \cite{PPZ3}, we had solved the equation (\ref{FY-form-eq}) for general dimension $n$ when $\alpha'<0$. However, the case $\alpha'>0$ for general dimension $n$ remained open, as a key lower bound for the Hessian could not be established \cite{PPZ1}. In this paper, we shall simultaneously solve the open case $\alpha'>0$ for general dimension $n$, improve on the solution found in \cite{PPZ3} for the case $\alpha'<0$, and do it actually for more general equations where the factor $(i\ddb u)^2$ in (\ref{FY-form-eq}) is replaced by higher powers of $i\ddb u$.

\smallskip
More precisely, let $(X, \hat\o)$, $\rho$, $\mu$, $\alpha'$ be as above. For each fixed integer $k$, $1\leq k \leq n-1$ and each real number $\gamma>0$, we consider the equation
\be \label{gen-eq-forms-intro}
i \ddb \left\{ e^{ku} \hat{\omega} - \alpha' e^{(k-\gamma)u} \rho \right\} \wedge \hat{\omega}^{n-2} + \alpha' (i \ddb u)^{k+1} \wedge \hat{\omega}^{n-k-1} + \mu \, \hat{\omega}^n = 0.
\ee
Clearly, when $k=1$ and $\gamma=2$, this equation reduces to the Fu-Yau equation (\ref{FY-form-eq}). We shall refer to (\ref{gen-eq-forms-intro}) as Fu-Yau Hessian equations. Our main result is then the following:

\begin{theorem} \label{main-thm}
Let $\alpha' \in {\bf R}$, $\rho \in \Omega^{1,1}(X,{\bf R})$, and $\mu: X \rightarrow {\bf R}$ be a smooth function such that $\int_X \mu \, \hat{\omega}^n = 0$. Define the set $\Upsilon_k$ by
\be \label{large-radius-regime}
\Upsilon_k = \left\{ u \in C^2(X,{\bf R}) : e^{-\gamma u} < \delta, \ |\alpha'| |e^{-u} i \ddb u|^k_{\hat{\omega}} < \tau \right\},
\ee
where $0<\delta, \, \tau \ll 1$ are explicit fixed constants depending only on $(X,\hat{\omega}), \alpha', \rho, \mu, n, k,\gamma$, whose expressions are given in (\ref{Gamma-cond}, \ref{delta-def}) below. Then 
there exists $M_0 \gg 1$ depending on $(X,\hat{\omega})$, $\alpha'$, $n$, $k$, $\gamma$, $\mu$ and $\rho$, such that for each $M \geq M_0$, there exists a unique smooth function $u \in \Upsilon_k$ with normalization $\int_X e^u \, \hat{\omega}^n = M$ solving the Fu-Yau Hessian equation {\rm (\ref{gen-eq-forms-intro})}.
\end{theorem}

We outline now the key differences between the earlier approaches 
and the approach of the present paper. 

\smallskip
The earlier approaches \cite{FY1, FY2, PPZ1, PPZ2} were based on rewriting the equation (\ref{FY-form-eq}) as
\bea
\label{FY-form-eq1}
\hat \sigma_2(\omega')={n(n-1)\over 2}
(e^{2u}-4\alpha'e^u|\na u|^2)+\nu
\eea
where $\nu$ is a linear combination of known functions, $u$ and $\na u$, $\omega'$ is defined by $\omega'=e^u\hat\omega+\alpha'e^{-u}\rho+2n\alpha'i\ddb u$, and $\hat\sigma_k(\o')$ is the $k$-th symmetric function of the eigenvalues of $\o'$ with respect to $\hat\o$. We look then for solutions $u$ satisfying the condition 
$\omega'\in\Gamma_2$, where $\Gamma_2$ is defined by the conditions $\hat\sigma_1(\omega')>0$ and $\hat\sigma_2(\o')>0$. The left hand side is then $>0$. When $\alpha'>0$, this implies immediately an upper bound on $|\na u|$. However, the difficulty is then to derive a positive lower bound for $\hat\sigma_2(\omega')$, and the arguments of \cite{FY1} worked only when $n=2$. On the other hand, when $\alpha'<0$, such a lower bound turns out to hold because there is no cancellation in the expression $e^{2u}-4\alpha'e^u|\na u|^2$. The estimate for $|\na u|$
and $|\hat\sigma_2(\o')|$ can then be obtained respectively by applying the techniques of Dinew-Kolodziej \cite{DK}, and Chou-X.J. Wang \cite{CW}, Hou-Ma-Wu \cite{HMW}, Guan \cite{Guan}, and the authors \cite{PPZ4}. 

\smallskip
The approach in the present paper relies instead on a different strategy. 

\smallskip

First, the equation (\ref{FY-form-eq}) corresponds to the case $k=1$, $\gamma=2$ of the Fu-Yau Hessian equations. As stated in Theorem \ref{main-thm}, we look for solutions $u\in \Upsilon_1$,
which is a more stringent condition than $\omega'\in\Gamma_2$. The set $\Upsilon_1$ and its condition $e^{-u} |\alpha' i \ddb u|_{\hat{\omega}} < \tau $ are inspired by the condition $|\alpha'Rm(\o)|<<1$ in \cite{PPZ3,PPZ4} which guarantees the parabolicity of the geometric flows introduced in these papers
\footnote{In these flows, a Hermitian metric $\o$ evolves with time, and $Rm(\o)$ is the curvature of the Chern unitary connection of $\o$. The condition $|\alpha'Rm(\o)|<<1$ was subsequently also used in \cite{FHP2}.)}. In the method of continuity, 
the given equation (\ref{FY-form-eq}) is realized as the end point of a family of equations for each $t\in [0,1]$. The condition $u\in\Upsilon_1$ implies that the diffusion operator $F^{p\bar q}\na_p\na_{\bar q}$ governing the evolution of $|Du|^2$ and $|\alpha'i\ddb u|^2$ is a controllable perturbation of the Laplacian $\Delta=g^{p\bar q}\na_p\na_{\bar q}$. The main problem is then to show that, if $u\in\Upsilon_1$ at time $t=0$, it will stay in $\Upsilon_1$ at all times.

\smallskip

This is accomplished by establishing a priori estimates, which we shall refer to as ``estimates with scale", which are more precise and delicate than the usual ones. Indeed, a priori estimates for $|u|$, $|Du|$, and $|\alpha'i\ddb u|$ are usually required only to be independent of $z\in X$ and $t\in[0,1]$. In the present situation, the normalization as 
given in Theorem \ref{main-thm}
\bea
\int_X e^u\hat\o^n=M
\eea
sets effectively a scale $M$, and the estimates with scale that we need are estimates for $|u|$, $|Du|$, and $|\alpha' i\ddb u|$ in terms of some specific powers of $M$. An example of such an estimate is the $C^0$ estimate stated in Theorem \ref{c0-est} below, $C^{-1}M\leq e^u\leq C\,M$,
which is a version in the present context of similar $C^0$ estimates established earlier in \cite{FY1, FY2, PPZ2}. The hardest part of the paper resides in the proof of similar estimates with scale for $|Du|$ and $|i\ddb u|$, as stated in Theorems \ref{c1-est} and \ref{c2-est}. 
Neither the set $\Upsilon_1$ nor the estimates with scale depend on the sign of $\alpha'$, which is why both cases $\alpha'>0$ and $\alpha'<0$ can be treated simultaneously. Furthermore we obtain a solution $u\in\Upsilon_1$, which is better than a solution in $\Gamma_2$. A vital clue that a strategy based on $\Upsilon_1$ and estimates with scale could work was provided by the authors' earlier alternative proof \cite{PPZ3, PPZ4} by flow methods of the Fu-Yau theorem \cite{FY1, FY2} in dimension $n=2$. 

\smallskip
The power of the new method is even more evident when it comes to the general Fu-Yau Hessian equation (\ref{gen-eq-forms-intro}). For $k\geq 2$, it is no longer possible to express the equation (\ref{gen-eq-forms-intro}) in terms of a single Hessian $\hat\sigma_{k+1}(\o')$ for some $(1,1)$-form $\o'$ as in (\ref{FY-form-eq1}). Rather, the equation leads to a combination of several Hessians, which makes it non-concave, and prevents a derivation of $C^2$ and $C^{2,\alpha}$ estimates by standard techniques of concave PDE's. On the other hand, the method of an ellipticity condition $\Upsilon_k$ preserved by estimates with scale works seamlessly in all cases of $1\leq k\leq n-1$. In fact the $C^3$ estimates that we obtain appear to be the first $C^3$ estimates established in the literature for any general class of Hessian equations besides the Laplacian and the Monge-Amp\`ere equations.

\section{Proof of Theorem 1: A Priori Estimates}
\setcounter{equation}{0}

In our study of  (\ref{gen-eq-forms-intro}), we will assume that ${\rm Vol}(X,\hat{\omega})=1$, which can be acheived by scaling $\hat{\omega} \mapsto \lambda \hat{\omega}$, $\alpha' \mapsto \lambda^k \alpha'$, $\rho \mapsto \lambda^{-k+1} \rho$, $\mu \mapsto \lambda^{-1} \mu$. Since the equation (\ref{gen-eq-forms-intro}) reduces to the Laplace equation when $\alpha'=0$, we assume from now on that $\alpha'\not=0$. We will use the notation $C^\ell_n = {n! \over \ell! (n-\ell)!}$ and $\hat{\sigma}_\ell(i \ddb u) \, \hat{\omega}^n = C^\ell_n \, (i \ddb u)^\ell \wedge \hat{\omega}^{n-\ell}$. Given $\rho$, we define the differential operator $L_\rho$ acting on functions by 
\bea
L_\rho f \, \hat{\omega}^n = n i \ddb (f \rho) \wedge \hat{\omega}^{n-2}. 
\eea
For each fixed $k \in \{1,2,3, \dots, n-1 \}$ and a real number $\gamma > 0$, the Fu-Yau Hessian equation (\ref{gen-eq-forms-intro}) can be rewritten as
\be \label{FY-scalar-eq}
{1 \over k} \, \Delta_{\hat{g}} e^{ku} + \alpha' \bigg\{ L_\rho e^{(k-\gamma)u} + \hat{\sigma}_{k+1}(i \ddb u) \bigg\} = \mu .
\ee
We note that we adjusted our conventions compared to the introduction by redefining $\mu$, $\rho$, and $\alpha'$ up to a constant. From this point on, we only work with the present conventions (\ref{FY-scalar-eq}). The standard Fu-Yau equation can be recovered by letting $k=1$, $\gamma = 2$. We remark that this equation is already of interest in the case when $\rho \equiv 0$, in which case the term $L_\rho e^{(k-\gamma)u}$ vanishes.
\smallskip
\par We can also write $L_\rho$ as
\be
L_\rho = a^{j \bar{k}} \p_j \p_{\bar{k}} + b^i \p_i + b^{\bar{i}} \p_{\bar{i}} + c,
\ee
where $a^{j \bar{k}}$ is a Hermitian section of $(T^{1,0}X)^* \otimes (T^{0,1}X)^*$, $b^i$ is a section of $(T^{1,0}X)^*$, and $c$ is a real function. All these coefficients are characterized by the following equations
\be
n i \ddb f \wedge \rho \wedge \hat{\omega}^{n-2} = a^{j \bar{k}} \p_j \p_{\bar{k}} f \, \hat{\omega}^n, \ \ n i \p f \wedge \bar{\p} \rho \wedge \hat{\omega}^{n-2} = b^i \p_i f \, \hat{\omega}^n, \ \ n i \ddb \rho \wedge \hat{\omega}^{n-2} = c \hat{\omega}^n.
\ee
for an arbitrary function $f$, and can be expressed explicitly in terms of $\rho$ and $\hat\o$ if desired. We will use the constant $\Lambda$ depending on $\rho$ defined by
\be
-\Lambda \hat{g}^{j \bar{k}} \leq a^{j \bar{k}} \leq \Lambda \hat{g}^{j \bar{k}}, \ \ \hat{\omega} =  \hat{g}_{\bar{k} j} i dz^j \wedge d \bar{z}^k, \ \ \hat{g}^{j \bar{k}} = (\hat{g}_{\bar{k} j})^{-1}.
\ee
We will look for solutions in the region
\be \label{Gamma-cond}
\Upsilon_k = \left\{ u \in C^2(X,{\bf R}) : e^{-\gamma u} < \delta, \ |\alpha'| |e^{-u} i \ddb u|^k_{\hat{\omega}} < \tau \right\}, \ \ \tau=  {2^{-7} \over C^k_{n-1}},
\ee
where $0<\delta \ll 1$ is a fixed small constant depending only on $(X,\hat{\omega}), \alpha', \rho, \mu, k, n, \gamma$. More precisely, it suffices for $\delta$ to satisfy the inequality
\be \label{delta-def}
\delta \leq \min\left\{1, {2^{-13} \over  |\alpha'| (k+\gamma)^3 \Lambda}, \bigg( {\theta \over 2 C_X \left(\| \mu \|_\infty+ \|\alpha' c \|_\infty \right)}  \bigg)^{\gamma / \gamma'} \right\},
\ee
where
\be \label{C_1-defn}
\theta = {1 \over 2C_1 -1}, \ \ \gamma' = \min \{ k, \gamma \}, \ \ C_1 = \{ 2(C_X+1)(\gamma+k)\}^{n} \bigg( {n \over n-1} \bigg)^{n^2}.
\ee
Here $C_X$ is the maximum of the constants appearing in the Poincar\'e inequality and Sobolev inequality on $(X,\hat{\omega})$. The proof of Theorem \ref{main-thm} is based on the following a priori estimates:

\begin{theorem}
\label{apriori}
Let $u \in \Upsilon_k$ be a $C^{5,\beta}(X)$ function with normalization $\int_X e^u \, \hat{\omega}^n = M$ solving the $k$-th Fu-Yau Hessian equation {\rm (\ref{FY-scalar-eq})}. Then
\be \label{apriori-estimates}
C^{-1} M \leq e^u \leq C M, \ \ e^{-u} | i \ddb u|_{\hat{\omega}} \leq C M^{-1/2}, \ \ e^{-3u} | \na \bar{\na} \na u |^2_{\hat{\omega}} \leq C',
\ee
where $C>1$ only depends on $(X,\hat{\omega})$, $\alpha'$, $k$, $\gamma$, $n$, $\rho$, and $\mu$. 
\end{theorem}

\medskip
Assuming Theorem \ref{apriori}, we can prove
Theorem \ref{main-thm}. Both the existence and uniqueness statements will be proved by the continuity method. We begin with the existence. Fix $\alpha' \in {\bf R} \backslash \{ 0 \}$, $\gamma>0$, $1 \leq k \leq (n-1)$, $\rho \in \Omega^{1,1}(X,{\bf R})$ and $\mu: X \rightarrow {\bf R}$ such that $\int_X \mu \, \hat{\omega}^n=0$, and define the set $\Upsilon_k$ as above. For a real parameter $t$, we consider the family of equations
\be \label{cont-method-scalar}
{1 \over k} \, \Delta_{\hat{g}} e^{ku_t} + \alpha' \left\{ t L_\rho e^{(k-\gamma)u_t} + \hat{\sigma}_{k+1}(i \ddb u_t) \right\} = t \mu .
\ee
As equations of differential forms, this family can be expressed as
\be \label{cont-method-form}
i \ddb \left\{ {e^{ku} \over k} \hat{\omega} + \alpha' t e^{(k-\gamma)u} \rho \right\} \wedge \hat{\omega}^{n-2} +  \alpha'   {C^{k}_{n-1} \over k+1} (i \ddb u)^{k+1} \wedge \hat{\omega}^{n-k-1} - t {\mu \over n} \hat{\omega}^n = 0.
\ee
We introduce the following spaces
\be
B_M = \{ u \in C^{5,\beta}(X,{\bf R}) : \int_X e^u \, \hat{\omega}^n = M \},
\ee
\be
B_1 = \{ (t,u) \in [0,1] \times B_M : u \in \Upsilon_k \},
\ee
\be
B_2 = \{ \psi \in C^{3,\beta}(X,{\bf R}) : \int_X \psi \, \hat{\omega}^n  = 0 \}
\ee
and define the map $\Psi: B_1 \rightarrow B_2$ by
\be
\Psi(t,u) = {1 \over k} \, \Delta_{\hat{g}} e^{ku_t} + \alpha' t L_\rho e^{(k-\gamma)u_t} + \alpha' \hat{\sigma}_{k+1}(i \ddb u_t) - t \mu.
\ee
We consider
\be
I = \{ t \in [0,1] : {\rm there \ exists \ } u \in B_M \ {\rm such \ that} \ (t,u) \in B_1 \ {\rm and} \ \Psi(t,u) =0 \}. 
\ee

\smallskip
First, $0\in I$: indeed the constant function $u_0 = \log M - \log \int_X \hat{\omega}^n$ is in $\Upsilon_k$ when $M \gg 1$, and $u_0$ solves the equation at $t=0$. In particular $I$ is non-empty. 

\smallskip
Next, we show that $I$ is open.
Let $(t_0,u_0) \in B_1$, and let $L= (D_u \Psi)_{(t_0,u_0)}$ be the linearized operator at $(t_0,u_0)$,
\be
L : \bigg\{ h \in C^{5,\beta}(X,{\bf R}) : \int_X h e^{u_0} \, \hat{\omega}^n = 0 \bigg\} \rightarrow \bigg\{ \psi \in C^{3,\beta}(X,{\bf R}) : \int_X \psi \, \hat{\omega}^n = 0 \bigg\},
\ee
defined by
\bea
L(h) \hat{\omega}^n &=& i \ddb \{ e^{k u_0} h \, \hat{\omega} + \alpha' (k-\gamma) t_0 e^{(k-\gamma)u_0}h \, \rho \} \wedge \hat{\omega}^{n-2} \nonumber\\
&&+ \alpha' C^k_{n-1} i \ddb h \wedge (i \ddb u_0)^k \wedge \hat{\omega}^{n-k-1}.
\eea
The leading order terms are
\be
 L(h) \hat{\omega}^n =  e^{k u_0} \chi_{(t_0,u_0)} \wedge \hat{\omega}^{n-k-1} \wedge i \ddb h + \cdots
\ee
where 
\be
\chi_{(t,u)} = \hat{\omega}^k + \alpha' (k-\gamma) t e^{-\gamma u} \, \rho \wedge \hat{\omega}^{k-1} + \alpha' C^k_{n-1} (e^{-u} i \ddb u)^k.
\ee
Since $u_0 \in \Upsilon_k$, we see from the conditions (\ref{Gamma-cond}) that $\chi_{(t_0,u_0)} > 0$ as a $(k,k)$ form and hence $L$ is elliptic. The $L^2$ adjoint $L^*$ is readily computed by integrating by parts:
\bea
\int_X \psi L(h) \, \hat{\omega}^n &=& \int_X h \,  e^{ku_0} \chi_{(t_0,u_0)} \wedge \hat{\omega}^{n-k-1} \wedge  i \ddb \psi \nonumber\\
&=& \int_X h L^*(\psi) \, \hat{\omega}^n.
\eea
Since $L^*$ is an elliptic operator with no zeroth order terms, by the strong maximum principle the kernel of $L^*$ consists of constant functions. An index theory argument (see e.g. \cite{PPZ3} or \cite{FY1} for full details) shows that the kernel of $L$ is spanned by a function of constant sign. It follows that $L$ is an isomorphism. By the implicit function theorem, there exists a unique solution $(t,u_t)$ for $t$ sufficiently close to $t_0$, with $u_t\in\Upsilon_k$ since $\Upsilon_k$ is open. We conclude that $I$ is open.

\smallskip

Finally, we apply Theorem \ref{apriori} to show that $I$ is closed. Consider a sequence $t_i \in I$ such that $t_i \rightarrow t_\infty$, and denote $u_{t_i} \in \Upsilon_k \cap B_M$ the associated $C^{5,\beta}$ functions such that $\Psi(t_i,u_{t_i})=0$. By differentating the equation $e^{-k u_{t_i}} \Psi(t_i,u_{t_i})=0$ with the Chern connection $\hat{\na}$ of the K\"ahler metric $\hat{\omega}$, we obtain
\bea \label{differentiated-eqn}
0 &=& {\chi_{(t_i,u_{t_i})} \wedge \hat{\omega}^{n-k-1} \wedge i \ddb (\p_\ell u_{t_i}) \over \hat{\omega}^n} \nonumber\\
&&+\hat{\na}_\ell \{\alpha' t_i e^{-\gamma u_{t_i}}((k-\gamma)^2 a^{p \bar{q}} \p_p u_{t_i} \, \p_{\bar{q}} u_{t_i}  + (k-\gamma) b^k \p_k u_{t_i} + (k-\gamma) b^{\bar{k}} \p_{\bar{k}} u_{t_i} + c)\} \nonumber\\
&&+ k \hat{\na}_\ell |Du_{t_i}|^2_{\hat{g}}  - \alpha' k e^{-k u_{t_i}} \hat{\sigma}_{k+1}(i \ddb u_{t_i}) \p_\ell u_{t_i} - t_i \p_\ell \{e^{-ku_{t_i}} \mu \}.
\eea
Since the equations (\ref{cont-method-scalar}) are of the form (\ref{FY-scalar-eq}) with uniformly bounded coefficients $\rho$ and $\mu$, Theorem \ref{apriori} applies to give uniform control of $| u_{t_i}|$ and $| \p \bar{\p} \p u_{t_i}|_{\hat{\omega}}$ along this sequence. Therefore $\hat{\Delta} u_{t_i}$ is uniformly controlled in $C^\beta(X)$ for any $0<\beta<1$. By Schauder estimates, we have $\| u_{t_i} \|_{C^{2,\beta}} \leq C$. 

Thus the differentiated equation (\ref{differentiated-eqn}) is a linear elliptic equation for $\p_\ell u_{t_i}$ with $C^{\beta}$ coefficients. This equation is uniformly elliptic along the sequence, since $\chi_{(t_i,u_{t_i})} \geq {1 \over 2} \hat{\omega}^k$ by (\ref{apriori-estimates}) when $M \gg 1$. By Schauder estimates, we have uniform control of $\| D u_{t_i} \|_{C^{2,\beta}}$. A bootstrap argument shows that we have uniform control of $\| u_{t_i} \|_{C^{6,\beta}}$, hence we may extract a subsequence converging to $u_\infty \in C^{5,\beta}$. Furthermore, for $M\geq M_0 \gg 1$ large enough, we see from (\ref{apriori-estimates}) that
\be
e^{-u_\infty} \ll 1, \ \ |e^{-u} i \ddb u_\infty|_{\hat{\omega}} \ll 1,
\ee
hence $u_\infty \in \Upsilon_k$. Thus $I$ is closed.

\smallskip
\par Hence $I=[0,1]$ and consequently there exists a $C^{5,\beta}$ function $u \in \Upsilon_k$ with normalization $\int_X e^u \, \hat{\omega}^n=M$ solving the Fu-Yau equation (\ref{FY-scalar-eq}). By applying Schauder estimates and a bootstrap argument to the differentiated equation (\ref{differentiated-eqn}), we see that $u$ is smooth. 

\smallskip
We complete now the proof of Theorem \ref{main-thm} with the proof of uniqueness.
\smallskip
\par First, we show that the only solutions of the equation 
\be
{1 \over k} \, i \ddb e^{ku} \wedge \hat{\omega}^{n-1} + \alpha' {C^k_{n-1} \over k+1} (i \ddb u)^{k+1} \wedge \hat{\omega}^{n-k-1} =0
\ee
with $|\alpha'| C^k_{n-1} |e^{-u} i \ddb u |^k_{\hat{\omega}} <2^{-7}$ are constant functions. Multiplying by $u$ and integrating, we see that
\be
0 = \int_X i \p u \wedge \bar{\p} u \wedge \bigg\{ e^{ku} \hat{\omega}^k + \alpha' {C^k_{n-1} \over k+1} (i \ddb u)^k \bigg\} \wedge \hat{\omega}^{n-k-1},
\ee
and hence $u$ must be constant since $e^{ku} \hat{\omega}^k + \alpha' {C^k_{n-1} \over k+1} (i \ddb u)^k >0$ as a $(k,k)$ form.
\smallskip
\par Now suppose there are two distinct solutions $u \in \Upsilon_k$ and $v \in \Upsilon_k$ satisfying (\ref{FY-scalar-eq}) under the normalization $\int_X e^u \, \hat{\omega}^n = \int_X e^v \, \hat{\omega}^n = M$ with $M \geq M_0$. For $t \in [0,1]$, define
\bea
\Phi(t,u) &=& i \ddb \left\{ {e^{ku} \over k} \hat{\omega} + \alpha' (1-t) e^{(k-\gamma)u} \rho \right\} \wedge \hat{\omega}^{n-2} \nonumber\\
&&+  \alpha'   {C^{k}_{n-1} \over k+1} (i \ddb u)^{k+1} \wedge \hat{\omega}^{n-k-1} - (1-t) {\mu \over n} \hat{\omega}^n,
\eea 
and consider the path $t \mapsto u_t$ satisfying $\Phi(t,u_t)=0$, $u_t \in \Upsilon_k$, $\int_X e^{u_t} \hat{\omega}^n=M$ with initial condition $u_0=u$. 
\smallskip
\par The same argument which shows that $I$ is open also shows that the path $u_t$ exists for a short-time: there exists $\epsilon>0$ such that $u_t$ is defined on $[0,\epsilon)$. By our estimates (\ref{apriori-estimates}), we may extend the path to be defined for $t \in [0,1]$. By uniqueness of the equation with $t=1$, we know that $u_1 = \log M - \log \int_X \hat{\omega}^n$. The same argument gives a path $t \mapsto v_t$ satisfying $\Phi(t,v_t)=0$, $v_t \in \Upsilon_k$, $\int_X e^{v_t} \hat{\omega}^n=M$ with $v_0=v$ and $v_1 = \log M - \log \int_X \hat{\omega}^n$. But then at the first time $0<t_0 \leq 1$ when $u_{t_0}=v_{t_0}$, we contradict the local uniqueness of $\Phi(t,u_t)=0$ given by the implicit function theorem.
\smallskip
\par It follows from our discussion that in order to prove Theorem \ref{main-thm}, it remains to establish the a priori estimates (\ref{apriori-estimates}).

\section{The Uniform Estimate}
\setcounter{equation}{0}

\begin{theorem}\label{c0-est}
Suppose $u \in \Upsilon_k$ solves (\ref{FY-scalar-eq}) subject to the normalization $\int_X e^u \, \hat{\omega}^n = M$.  Then
\be
C^{-1} M \leq e^u \leq CM,
\ee
where $C$ only depends on $(X,\hat{\omega})$, $k$, and $\gamma$. 
\end{theorem}

We first note the following general identity which holds for any function $u$.
\be \label{gen-eq-forms}
0 = \alpha' (p-k) \int_X e^{(p-k)u} i\p u \wedge \bar\p u \wedge (i \ddb u)^{k} \wedge \hat{\omega}^{n-k-1}  + \alpha' \int_X e^{(p-k)u} \, (i\ddb u)^{k+1} \wedge \hat\omega^{n-k-1}.
\ee
Substituting the Fu-Yau Hessian equation (\ref{cont-method-form}) with $t=1$, we obtain
\bea
0&=& \alpha' {C^{k}_{n-1} \over k+1} (p-k) \int_X e^{(p-k) u} i\p u \wedge \bar\p u \wedge (i \ddb u)^{k} \wedge \hat{\omega}^{n-k-1}  \nonumber\\
&& + \int_X e^{(p-k)u} {\mu \over n} - \int_X e^{(p-k)u} i \ddb \left\{ {e^{ku} \over k}\hat{\omega} + \alpha' e^{(k-\gamma)u} \rho \right\} \wedge \hat{\omega}^{n-2} .
\eea
We integrate by parts to derive
\bea
0&=& \alpha' {C^{k}_{n-1} \over k+1 }(p-k) \int_X e^{(p-k) u} i\p u \wedge \bar\p u \wedge (i \ddb u)^{k} \wedge \hat{\omega}^{n-k-1}  \nonumber\\
&& + \int_X e^{(p-k)u} {\mu \over n} + (p-k) \int_X  e^{pu} \, i \p u \wedge i \bar{\p} u  \wedge \hat{\omega}^{n-1} \nonumber\\
&&+ (p-k) \alpha' \int_X  e^{(p-k)u} \, i \p u \wedge i \bar{\p} ( e^{(k-\gamma)u} \rho ) \wedge \hat{\omega}^{n-2}.
\eea
Integrating by parts again gives
\bea \label{c0-identity}
& \ & (p-k) \int_X e^{pu}  i \p u \wedge \bar\p u \wedge \hat\omega^{n-k-1} \wedge \chi \nonumber\\
 &=& -\int_X e^{(p-k)u} \mu {\hat{\omega} \over n}  + {p-k  \over p-\gamma} \alpha' \int_X  e^{(p-\gamma)u} \wedge i \ddb  \rho \wedge \hat{\omega}^{n-2},
\eea
where we now assume $p > \gamma$ and we define
\be
\chi = \hat{\omega}^{k} + \alpha'(k-\gamma) e^{-\gamma u} \rho \wedge \hat{\omega}^{k-1} +  \alpha' {C^{k}_{n-1} \over k+1} (e^{-u} i \ddb u)^{k}.
\ee
Next, we estimate
\bea
i \p u \wedge \bar\p u \wedge \hat\omega^{n-k-1} \wedge \chi &=& {|\na u|^2_{\hat{\omega}} \over n} \, \hat{\omega}^{n} + \alpha'(k-\gamma) e^{-\gamma u} {a^{i \bar{j}} u_i u_{\bar{j}} \over n} \, \hat{\omega}^{n} \nonumber\\
&& +  \alpha' {C^{k}_{n-1} \over k+1}  i \p u \wedge \bar\p u \wedge  (e^{-u} i \ddb u)^{k} \wedge \hat\omega^{n-k-1} \nonumber\\
&\geq& {|\na u|^2_{\hat{\omega}} \over n} \, \hat{\omega}^{n} - |\alpha' \Lambda (k-\gamma)| \delta {|\na u|^2_{\hat{\omega}} \over n} \, \hat{\omega}^{n} \nonumber\\
&& - |\alpha'| {C^{k}_{n-1} \over k+1 } |e^{-u} i \ddb u|_{\hat{\omega}}^{k} {|\na u|^2_{\hat{\omega}} \over n} \, \hat\omega^{n}.
\eea
Since $u \in \Upsilon_k$, by (\ref{Gamma-cond}) and (\ref{delta-def}) the positive term dominates the expression and we can conclude
\be \label{gamma-c0-cond}
i \p u \wedge \bar\p u \wedge \hat\omega^{n-k-1} \wedge \chi \geq {1 \over 2} {|\na u|^2_{\hat{\omega}} \over n} \, \hat{\omega}^{n}.
\ee
The proof of Theorem \ref{c0-est} will be divided into three propositions. We note that in the following arguments we will omit the background volume form $\hat{\omega}^n$ when integrating scalar functions.
\begin{proposition} \label{e^u-leq-CM}
Suppose $u \in \Upsilon_k$ solves (\ref{FY-scalar-eq}) subject to normalization $\int_X e^u = M$. There exists $C_1>0$ such that
\be
e^u \leq C_1 M,
\ee
where $C_1$ only depends on $(X,\hat{\omega})$, $n$, $k$ and $\gamma$. In fact, $C_1$ is given by (\ref{C_1-defn}).
\end{proposition}

Combining (\ref{c0-identity}) and (\ref{gamma-c0-cond}) gives
\bea
& \ & {1 \over 2} (p-k) \int_X e^{pu} |\na u|^2_{\hat{\omega}} \nonumber\\
 &\leq& -\int_X e^{(p-k)u} \mu + {p-k  \over p-\gamma} n \alpha' \int_X  e^{(p-\gamma)u} \wedge i \ddb  \rho \wedge \hat{\omega}^{n-2}.
\eea
We estimate
\be
\int_X |\na e^{{p \over 2} u} |^2_{\hat{\omega}} \leq {p^2 \over 2 (p-k)} \bigg\{ \| \mu \|_{L^\infty} \int_X e^{(p-k)u} + {p-k  \over p-\gamma} \| \alpha' c\|_{L^\infty} \int_X  e^{(p-\gamma)u} \bigg\}.
\ee
For any $p \geq 2\max \{ \gamma,k \}$, there holds ${p^2 \over 2 (p-k)} \leq p$ and ${p-k \over p-\gamma} \leq 2$. Using $e^{-\gamma u} \leq \delta \leq 1$ and (\ref{delta-def}), we conclude that
\bea
\int_X |\na e^{{p \over 2} u} |^2_{\hat{\omega}} &\leq& 2(\| \mu \|_{L^\infty}+  \| \alpha' c\|_{L^\infty}) \delta^{\min \{ k,\gamma \} \over \gamma} p \int_X e^{pu} \nonumber\\
&\leq& {\theta \over C_X} p \int_X e^{pu} \leq {p \over C_X}  \int_X e^{pu},
\eea
for any $p \geq 2 (\gamma + k)$. Let $\beta = {n \over n-1}$. The Sobolev inequality gives us
\be
\left( \int_X e^{\beta pu} \right)^{1/\beta} \leq C_X \bigg( \int_X |\na e^{{p \over 2} u}|^2_{\hat{\omega}} + \int_X e^{pu} \bigg).
\ee
Therefore for all $p \geq 2 (\gamma+k)$,
\be
\|e^u\|_{L^{p\beta}} \leq (C_X+1)^{1/p} p^{1/ p} \|e^u \|_{L^p}.
\ee
Iterating this inequality gives
\bea
\|e^u\|_{L^{p\beta^{(k+1)}}} \leq \{ (C_X+1) p \}^{{1\over p} \sum_{i=0}^k{1\over \beta^{i}}} \cdot \beta^{{1\over p} \sum_{i=1}^k {i\over \beta^{i}}}\, \|e^u\|_{L^p}.
\eea
Letting $k\rightarrow \infty$, we obtain
\be
\sup_X e^u \leq C_1' \| e^u \|_{L^{2(\gamma +k)}}, \ \ C_1' = \{2 (C_X+1) (\gamma+k)\}^{{1\over 2 (\gamma+k)} \sum_{i=0}^\infty{1\over \beta^{i}}} \cdot \beta^{{1\over 2 (\gamma+k)} \sum_{i=1}^\infty{i \over \beta^{i}}}.
\ee
It follows that
\be
\sup_X e^u \leq C_1' (\sup_X e^u)^{1-(2(\gamma+k))^{-1}} \left( \int_X e^u \right)^{1/2(\gamma+k)},
\ee
and we conclude that
\be
\sup_X e^u \leq C_1 \int_X e^u, \ \ C_1 = (C_1')^{2(\gamma+k)}.
\ee
This proves the estimate. As it will be needed in the future, we note that the precise form of $C_1$ agrees with the definition given in (\ref{C_1-defn}).

\begin{proposition} \label{int-e-minus-u}
Suppose $u \in \Upsilon_k$ solves (\ref{FY-scalar-eq}) subject to normalization $\int_X e^u = M$. There exists a constant $C$ only depending on $(X,\hat{\omega})$, $n$, $k$ and $\gamma$ such that
\be 
\int_X e^{-u} \leq C M^{-1}.
\ee
\end{proposition}

\smallskip

Setting $p=-1$ in (\ref{c0-identity}) gives
\bea
&& (k+1)\int_X e^{-u} i\p u\wedge \bar\p u \wedge \hat\omega^{n-k-1} \wedge \chi \\
&=&  \int_X e^{-(1+k)u} \mu {\hat\omega^n\over n} - {1+k \over 1+\gamma} \int_X e^{-(1-\gamma) u} i\ddb \rho \wedge \hat\omega^{n-2}\nonumber\\
&\leq & {1 \over n} \|\mu\|_{L^\infty} \int_X e^{-(1+k)u} + {1+k \over (1+\gamma)n} \|\alpha' c \|_{L^\infty} \int_X e^{-(1+\gamma)u}.
\eea
Since $u\in \Upsilon_k$, we may use (\ref{gamma-c0-cond}) and $e^{- \gamma u} \leq \delta \leq 1$ to obtain
\be
\int_X e^{-u} |\na u|^2_{\hat{\omega}} \leq  2 \delta^{\min \{ k,\gamma \} \over \gamma} \left(\| \mu \|_{L^\infty}+ \|\alpha' c \|_{L^\infty} \right)  \int_X e^{-u}.
\ee
By the Poincar\'e inequality
\be
\int_X e^{-u} - \left( \int_X e^{-u/2} \right)^2 \leq C_X \int_X |\na e^{-u/2}|^2_{\hat{\omega}}.
\ee
After using the definition of $\delta$ (\ref{delta-def}), it follows that
\be \label{theta-int}
\int_X e^{-u} \leq {1\over 1-{\theta \over 4}} \left( \int_X e^{-u/2} \right)^2. 
\ee
Let $U = \{x \in X : e^{u} \geq {M \over 2} \}$. From Proposition \ref{e^u-leq-CM}, and using ${\rm Vol}(X,\hat{\omega})=1$,
\be
M = \int_X e^u \leq C_1 M |U|  + (1-|U|) {M \over 2}.
\ee
Hence $|U| \geq \theta >0$, where we recall that $\theta$ was defined in (\ref{delta-def}). Using $|U| \geq \theta$ and (\ref{theta-int}), it was shown in \cite{PPZ3} that the estimate
\be
\int_X e^{-u} \leq {1 \over 1 - {\theta \over 4}} \left( 1 + {2 \over \theta} \right) \left( {2 \over \theta^2} \right) M^{-1}
\ee
follows.
\begin{proposition}
Suppose $u \in \Upsilon_k$ solves (\ref{FY-scalar-eq}) subject to the normalization $\int_X e^u = M$. There exists $C$ such that
\be
\sup_X e^{-u} \leq C M^{-1},
\ee
where $C$ only depends on $(X,\hat{\omega})$, $n$, $k$ and $\gamma$.
\end{proposition}
Exchanging $p$ for $-p$ in (\ref{c0-identity}) and using (\ref{gamma-c0-cond}) gives
\bea
&& (p+k) \int_X e^{-p u}  i\p u \wedge \bar\p u\wedge \hat\omega^{n-1} \\
&\leq &  2 \int_X e^{-(p+k)u} \mu {\hat\omega^n\over n} - 2 \alpha' {p+k\over p+\gamma} \int_X e^{-(p+\gamma)u} i\ddb \rho\wedge \hat\omega^{n-2}.\nonumber
\eea
By using $e^{\gamma u} \leq \delta \leq 1$, we obtain
\be
\int_X |\na e^{-{p \over 2} u}|^2_{\hat{\omega}} \leq {p^2 \over 2 (p+k)} \delta^{\min \{ k, \gamma \} \over \gamma} (\| \mu \|_{L^\infty} + {p+k \over p + \gamma} \| \alpha' c \|_{L^\infty} )  \int_X e^{-pu}.
\ee
We may use (\ref{delta-def}) to obtain a constant $C$ depending on $(X,\hat{\omega})$, $n$, $k$, and $\gamma$ such that
\be
\int_X |\na e^{-{p \over 2} u}|^2_{\hat{\omega}} \leq C p \int_X e^{-pu}.
\ee
for any $p \geq 1$. Using the Sobolev inequality and iterating in a similar way to Proposition \ref{e^u-leq-CM}, we obtain
\be
\sup_X e^{-u} \leq C \| e^{-u} \|_{L^1}.
\ee
Applying Proposition \ref{int-e-minus-u} gives the desired estimate.

\section{Setup and Notation}
\setcounter{equation}{0}

\subsection{The formalism of evolving metrics}
We come now to the key steps of establishing the gradient and the $C^2$ estimates. It turns out that, for these steps, it is more natural to view the equation
(\ref{FY-scalar-eq}) as an equation for the unknown, non-K\"ahler, Hermitian form
\bea
\o=e^u\hat\o
\eea
and to carry out calculations with respect to the Chern unitary connection $\na$ of $\o$. As usual, we identify the metrics $\hat{g}$ and $g$ via $\hat{\omega} = \hat{g}_{\bar{k} j} \, i dz^j \wedge d\bar{z}^k$ and $\omega = g_{\bar{k} j} \, i dz^j \wedge d\bar{z}^k$, and denote $\hat{g}^{j \bar{k}}$, $g^{j \bar{k}}$ to be the inverse matrix of $\hat{g}_{\bar{k} j}$, $g_{\bar{k} j}$. Then $g_{\bar kj}=e^u\hat g_{\bar kj}$,
$g^{j\bar k}=e^{-u}\hat g^{j\bar k}$. Recall that the Chern unitary connection $\na$ is defined by
\bea
\na_{\bar k}V^j=\p_{\bar k}V^j,
\quad
\na_kV^j=g^{j\bar m}\p_k(g_{\bar m p}V^p)
\eea
and its torsion and curvature by
\bea \label{curv-torsion-defn}
[\na_\alpha,\na_\beta]V^\gamma=R_{\beta\alpha}{}^\gamma{}_\delta V^\delta
+
T^\delta{}_{\beta\alpha}\na_\delta V^\gamma.
\eea
Explicitly,
\bea
R_{\bar kq}{}^j{}_p=-\p_{\bar k}(g^{j\bar m}\p_jg_{\bar m q}),
\quad
T^j{}_{pq}=g^{j\bar m}(\p_pg_{\bar m q}-\p_qg_{\bar m p}).
\eea
The curvatures and torsions of the metrics $g_{\bar kj}$ and $\hat g_{\bar kj}$ are then related by
\bea \label{curv-torsion}
R_{\bar{k} j}{}^p{}_i = \hat{R}_{\bar{k} j}{}^p{}_i - u_{\bar{k} j} \delta^p{}_i,
\quad
T^\lambda{}_{kj} = u_k \delta^\lambda{}_j - u_j \delta^\lambda{}_k.
\eea
The following commutation formulas with either $3$ or $4$ covariant derivatives will be particularly useful,
\bea
\na_j\na_p\na_{\bar q}u
=
\na_p\na_{\bar q}\na_j u
+
T^m{}_{pj}\na_m\na_{\bar q}u
\eea
and
\bea
\na_{\bar k}\na_j\na_p\na_{\bar q}u
&=&
\na_p\na_{\bar q}\na_j\na_{\bar k}u-R_{\bar qp\bar k}{}^{\bar m}\na_{\bar m}\na_ju
+
R_{\bar kj}{}^m{}_p\na_m\na_{\bar q}u
\nonumber\\
&&
+T^{\bar m}{}_{\bar q\bar k}\na_p\na_{\bar m}\na_ju
+
T^m{}_{pj}\na_{\bar k}\na_m \na_{\bar q}u.
\eea
They reduce in our case to
\be \label{na-na-na-u}
\na_j \na_p \na_{\bar{q}} u =  \na_p \na_{\bar{q}} \na_j u + u_p u_{\bar{q} j}  - u_j u_{\bar{q} p}.
\ee
and to
\bea \label{na-na-na-na-u}
\na_{\bar{k}} \na_j \na_p \na_{\bar{q}} u &=& \na_p \na_{\bar{q}} \na_j \na_{\bar{k}} u + u_p \na_{\bar{k}} \na_j \na_{\bar{q}} u - u_j \na_{\bar{k}} \na_p \na_{\bar{q}} u \nonumber\\
&&+ u_{\bar{q}} \na_p  \na_{\bar{k}}  \na_j u -  u_{\bar{k}} \na_p \na_{\bar{q}}  \na_j u \nonumber\\
&& + \hat{R}_{\bar{k} j}{}^{\lambda}{}_p u_{\bar{q} \lambda}  - \hat{R}_{\bar{q} p \bar{k}}{}^{\bar{\lambda}} u_{\bar{\lambda} j}.
\eea

\smallskip

It will also be convenient to use the symmetric functions of the eigenvalues of $i\ddb u$ with respect to $\o$ rather than with respect to $\hat\o$. Thus we 
define $\sigma_\ell(i \ddb u)$ to be the $\ell$-th elementary symmetric polynomial of the eigenvalues of the endomorphism $h^{i}{}_{j}= g^{i \bar{k}} u_{\bar{k} j}$. Explicitly, if $\lambda_i$ are the eigenvalues of the endomorphism $h^{i}{}_{j} = g^{i \bar{k}} u_{\bar{k} j}$, then $\sigma_\ell(i \ddb u) = \sum_{i_1< \cdots < i_\ell} \lambda_{i_1} \cdots \lambda_{i_\ell}$. Using this formalism, equation (\ref{FY-scalar-eq}) becomes
\be \label{FY-scalar-g-compact}
\Delta_g u + k|Du|^2_g + \alpha' e^{-(k+1)u} L_\rho e^{(k-\gamma)u} + \alpha' \sigma_{k+1}(i \ddb u) - e^{-(k+1)u} \mu = 0.
\ee

\subsection{Differentiating Hessian operators}
We define
\be
\sigma_\ell^{p \bar{q}} =  {\p \sigma_\ell \over \p h^r{}_p} g^{r \bar{q}}, \ \ \sigma_\ell^{p \bar{q}, r \bar{s}} = {\p^2 \sigma_\ell \over \p h^a{}_p \p h^b{}_r} g^{a \bar{q}} g^{b \bar{s}}.
\ee
Then the variational formula $\delta \sigma_\ell =  {\p \sigma_\ell \over \p h^r{}_p} \delta h^r{}_p$ becomes
\be \label{diff-sigmak}
\na_i \sigma_\ell = \sigma_\ell^{p \bar{q}} \na_i u_{\bar{q} p}.
\ee
Similarly,
\be \label{diff-diff-sigmak}
\na_{\bar{j}} \sigma_\ell^{p \bar{q}} = \sigma_\ell^{p \bar{q}, r \bar{s}} \na_{\bar{j}} u_{\bar{s} r}. 
\ee
 We will use a general formula for differentiating a function of eigenvalues of a matrix. Let $F(h)= f(\lambda_1, \cdots, \lambda_n)$ be a symmetric function of the eigenvalues of a Hermitian matrix $h$. Then at a diagonal matrix $h$, we have (see \cite{Ball, Gerhardt}), 
\begin{eqnarray}
 {\p F \over \p h^i{}_j} &=& \delta_{ij} f_i,\\
\sum {\p^2 F \over \p h^i{}_j \p h^r{}_s} T^i{}_j T^r{}_{s} &=& \sum f_{ij} T^i{}_i T^j{}_j + \sum_{p\neq q}\frac{f_p - f_q}{\lambda_p-\lambda_q} | T^p{}_q|^2.
\end{eqnarray}
for any Hermitian matrix $T$. Since $\sigma_\ell(h)= \sum_{i_1 < \cdots < i_\ell} \lambda_{i_1} \lambda_{i_2} \cdots \lambda_{i_\ell}$, this formula implies that at a point $p \in X$ where $g$ is the identity and $u_{\bar{q} p}$ is diagonal, then
\begin{eqnarray}
 \label{firstorder} \sigma_\ell^{p \bar{q}} &=& \delta_{pq} \sigma_{\ell-1}(\lambda|p),\\
 \label{secondorder} \sigma_\ell^{p \bar{q}, r \bar{s}} \na_i u_{\bar{q} p} \na_{\bar{i}} u_{\bar{s} r} &=& \sum_{p,q} \sigma_{\ell-2}(\lambda|pq) \na_i u_{\bar{p} p} \na_{\bar{i}} u_{\bar{q} q} - \sum_{p\neq q} \sigma_{\ell-2}(\lambda|pq) | \na_i u_{\bar{q} p}|^2.
\end{eqnarray}
We introduced the notation $\sigma_m(\lambda|p)$ and $\sigma_m(\lambda|pq)$ for the $m$-th elementary symmetric polynomial of 
\bea
\nonumber(\lambda | i) = (\lambda_1, \cdots, \widehat{\lambda_{i}}, \cdots, \lambda_n)\in \R^{n-1} \ {\rm and } \
(\lambda | ij) = (\lambda_1, \cdots, \widehat{\lambda_{i}}, \cdots, \widehat{\lambda_{j}}, \cdots, \lambda_n)\in \R^{n-2}.
\eea
Lastly, we introduce the tensor $F^{p \bar{q}}$, which will appear in subsequent sections when we differentiate the Fu-Yau equation.
\be \label{defn-Fpq}
F^{p \bar{q}} =  g^{p \bar{q}} + \alpha' (k-\gamma) e^{-(1+\gamma)u} a^{p \bar{q}} +  \alpha' \sigma_{k+1}^{p \bar{q}}.
\ee
We will prove that for $u \in \Upsilon_k$, $F^{p \bar{q}}$ is close to the metric $g^{p \bar{q}}$. For this, we first note the following elementary estimate.

\begin{lemma} \label{sigma-est-lemma}
Let $m$ be a positive integer and $\ell \in \{1, \dots, m \}$. For any vector $\lambda \in {\bf R}^m$,
\be
|\sigma_\ell(\lambda) | \leq {C^\ell_m \over m^{\ell/2}} |\lambda|^\ell 
\ee
with $|\lambda|= \left(\sum_{i=1}^n \lambda_i^2\right)^{1/2}$. Here, $\sigma_\ell(\lambda)$ is the $\ell$-th elementary symmetric polynomial of $\lambda$ and $C^\ell_m = {m! \over \ell! (m-\ell)!}$. 
\end{lemma}

{\it Proof:} Using the Newton-Maclaurin inequality, 
\be
|\sigma_\ell(\lambda)| \leq \sigma_\ell(|\lambda_1|, \dots, |\lambda_m|) \leq C^\ell_m \left( {\sum_i^m |\lambda_i| \over m} \right)^\ell.
\ee
The Cauchy-Schwarz inequality now gives the desired estimate. Q.E.D.

\bigskip

\par We can now prove the following simple but important lemma regarding the ellipticity of $F^{p \bar{q}}$.

\begin{lemma}\label{F-est-lemma}
If $u\in \Upsilon_k$, then 
\be \label{F-2-g}
(1-2^{-6})g^{p \bar{q}} \leq F^{p \bar{q}} \leq (1+2^{-6})g^{p \bar{q}}.
\ee
\end{lemma}
\smallskip
{\it Proof:} 
First, at a point $z$ where $g^{p \bar{q}} = \delta_{pq}$ and $u_{\bar{q} p}$ is diagonal, the above lemma implies
\be
|\alpha' \sigma_{k+1}^{p \bar{p}}| = |\alpha' \sigma_k(\lambda|p)| \leq |\alpha'| {C^k_{n-1} \over (n-1)^{k/2}} |\na \bar{\na} u|^k_g.
\ee
The condition $u \in \Upsilon_k$ gives $|\alpha' \sigma_{k+1}^{p \bar{p}}(z)| \leq 2^{-7}$. This argument shows that $\alpha' \sigma_{k+1}^{p \bar{q}}$ is on the order of $2^{-7} g^{p \bar{q}}$ in arbitrary coordinates.

Next, $u \in \Upsilon_k$ also implies that $|\alpha' (k-\gamma) e^{-\gamma u} \Lambda| \leq 2^{-7}$. Since $-\Lambda \hat{g}^{p \bar{q}} \leq a^{p \bar{q}} \leq \Lambda \hat{g}^{p \bar{q}}$, we can put everything together and obtain the estimate (\ref{F-2-g}). Q.E.D.

\section{Gradient Estimate}
\setcounter{equation}{0}

\begin{theorem}
\label{c1-est}
Let $u \in \Upsilon_k$ be a $C^{3}(X,{\bf R})$ function solving the Fu-Yau Hessian equation (\ref{FY-scalar-eq}). Then
\be
|\na u|^2_{\hat{g}} \leq C,
\ee
where $C$ only depends on $(X,\hat{\omega})$, $\alpha'$, $k$, $\gamma$, $\| \rho \|_{C^3(X, \hat\omega)}$ and $\| \mu \|_{C^1(X)}$. 
\end{theorem}

In view of Theorem \ref{c0-est}, this estimate is equivalent to
\be
|\na u|^2_g \leq C M^{-1},
\ee
where $C$ only depends on $(X,\hat{\omega})$, $\alpha'$, $k$, $\gamma$, $\| \rho \|_{C^3(X, \hat\omega)}$ and $\| \mu \|_{C^1(X)}$. 
We will prove this estimate by applying the maximum principle to
the following test function
\be
G =\log |\na u|_g^2 + (1+\sigma) u,
\ee
for a parameter $0<\sigma<1$. Though there is a range of values of $\sigma$ which makes the argument work, to be concrete we will take $\sigma = 2^{-7}$. 

\subsection{Estimating the leading terms}
Suppose $G$ attains a maximum at $p \in X$. Then
\be \label{critical-eqn}
0 = {\na |\na u|_g^2 \over |\na u|_g^2} + (1+\sigma) \na u.
\ee
We will compute the operator $F^{p \bar{q}} \na_p \na_{\bar{q}}$ acting on $G$ at $p$.
\be \label{grad-main-ineq-1}
F^{p \bar{q}} \na_p \na_{\bar{q}} G = {1 \over |\na u|^2_g} F^{p \bar{q}} \na_p \na_{\bar{q}} |\na u|^2_g - {1 \over |\na u|^4_g} F^{p \bar{q}} \na_p |\na u|^2_g \na_{\bar{q}} |\na u|^2_g + (1+\sigma) F^{p \bar{q}} u_{\bar{q} p}.
\ee
By direct computation
\bea
F^{p \bar{q}} \na_p \na_{\bar{q}} |\na u|^2_g &=& F^{p \bar{q}} g^{j \bar{i}} \na_p \na_{\bar{q}} \na_j u \na_{\bar{i}} u +  F^{p \bar{q}} g^{j \bar{i}} \na_j u  \na_p \na_{\bar{q}} \na_{\bar{i}} u \nonumber\\
&&+ |\na \bar{\na} u|^2_{Fg} + |\na \na u|^2_{Fg}.
\eea
where $|\na \na u|^2_{Fg} = F^{p \bar{q}} g^{j \bar{i}} \na_p \na_j u \na_{\bar{q}} \na_{\bar{i}} u$ and $|\na \bar{\na} u|^2_{Fg} = F^{p \bar{q}} g^{j \bar{i}} u_{\bar{q} j} u_{\bar{i} p}$. Commuting derivatives according to the relation
\be
[\na_j, \na_{\bar{\ell}}] u_{\bar{i}} = R_{\bar{\ell} j}{}_{\bar{i}}{}^{\bar{p}} u_{\bar{p}} = \hat{R}_{\bar{\ell} j}{}_{\bar{i}}{}^{\bar{p}} u_{\bar{p}} - u_{\bar{\ell} j} u_{\bar{i}},
\ee
we obtain
\be \label{grad-est-curv-appear}
F^{p \bar{q}} g^{j \bar{i}} \na_j u  \na_p \na_{\bar{q}} \na_{\bar{i}} u = \overline{F^{p \bar{q}} g^{j \bar{i}} \na_p  \na_{\bar{q}} \na_{j} u  \na_{\bar{i}} u } + F^{p \bar{q}} g^{j \bar{i}} u_j \hat{R}_{\bar{q} p}{}_{\bar{i}}{}^{\bar{\lambda}} u_{\bar{\lambda}} - F^{p \bar{q}} g^{j \bar{i}} u_{j} u_{\bar{q} p} u_{\bar{i}}.
\ee
Thus
\bea \label{DD|Du|^2-1}
F^{p \bar{q}} \na_p \na_{\bar{q}} |\na u|^2_g &=& 2 \Re \{ F^{p \bar{q}} g^{j \bar{i}} \na_p \na_{\bar{q}} \na_j u \na_{\bar{i}} u \} + F^{p \bar{q}} g^{j \bar{i}} u_j \hat{R}_{\bar{q} p}{}_{\bar{i}}{}^{\bar{\lambda}} u_{\bar{\lambda}} \nonumber\\
&& - F^{p \bar{q}} g^{j \bar{i}} u_{j} u_{\bar{q} p} u_{\bar{i}} + |\na \bar{\na} u|^2_{Fg} + |\na \na u|^2_{Fg}.
\eea
Next, we use the equation. Expanding $L_\rho = a^{p \bar{q}} \p_p \p_{\bar{q}} + b^i \p_i + \bar{b}^i \p_{\bar{i}} + c$, equation (\ref{FY-scalar-g-compact}) becomes
\bea \label{FY-scalar-g}
0 &=& \Delta_g u + \alpha' \left\{(k-\gamma) e^{-(1+\gamma)u} a^{p \bar{q}} u_{\bar{q} p} +  \sigma_{k+1}(i \ddb u) \right\} + k|\na u|^2_g \nonumber\\
&&+ \alpha' (k-\gamma)^2 e^{-(1+\gamma)u} a^{p \bar{q}} u_p u_{\bar{q}} +2 \alpha' (k-\gamma) e^{-(1+\gamma)u} \Re \{ b^i u_i \} \nonumber\\
&&+\alpha' e^{-(1+\gamma)u} c  - e^{-(k+1)u} \mu .
\eea
We covariantly differentiate equation (\ref{FY-scalar-g}), using (\ref{diff-sigmak}) to differentiate $\sigma_{k+1}$ and using the notation $F^{p \bar{q}}$ introduced in (\ref{defn-Fpq}). This leads to
\be \label{diff-once}
0 = F^{p \bar{q}} \na_j \na_p \na_{\bar{q}} u + k \na_j |\na u|^2_g + \E_j,
\ee
where 
\bea \label{E_j}
{\mathcal{E}}_j &=&  \alpha'(k-\gamma)e^{-(1+\gamma)u}  \bigg\{ -\gamma  a^{p \bar{q}} u_{\bar{q} p} u_j + \hat{\na}_j a^{p \bar{q}} u_{\bar{q} p} \bigg\} \nonumber\\
&& + \alpha' (k-\gamma)^2 e^{-(1+\gamma)u}  \bigg\{ -\gamma a^{p \bar{q}} u_p u_{\bar{q}} u_j  +  \hat{\na}_j a^{p \bar{q}} u_p u_{\bar{q}} +  a^{p \bar{q}} \na_j \na_p u u_{\bar{q}} + a^{p \bar{q}} u_p u_{\bar{q} j} \bigg\} \nonumber\\
&& +\alpha' (k-\gamma)  e^{-(1+\gamma)u} \bigg\{ -2(1+\gamma) \Re \{ b^i u_i \} u_j  + \hat{\na}_j b^i u_i   \nonumber\\
&& +   u_j b^i u_i + \p_j \bar{b^i} u_{\bar{i}} +   b^i \na_j \na_i u  + \bar{b^i} u_{\bar{i} j} \bigg\} \nonumber\\
&&-(1+\gamma) \alpha' e^{-(1+\gamma)u} c u_j +\alpha' e^{-(1+\gamma)u} \p_j c \nonumber\\
&& + (k+1)e^{-(k+1)u} \mu u_j - e^{-(k+1)u} \p_j \mu.
\eea
We used $\na_i W^j = \hat{\na}_i W^j + u_i W^j$ to replace $\nabla$ by $\hat\nabla$ in the above calculation. We will eventually see that the terms $\E_j$ play a minor role when $u \in \Upsilon_k$, and will only perturb the coefficients of the leading terms. Commuting covariant derivatives using (\ref{na-na-na-u}), we obtain
\be \label{diff-eqn-once}
F^{p \bar{q}} \na_p \na_{\bar{q}} \na_j u = - F^{p \bar{q}} u_p u_{\bar{q} j} + F^{p \bar{q}} u_j u_{\bar{q} p} - k \na_j |\na u|^2_g - \E_j.
\ee
Substituting (\ref{diff-eqn-once}) into (\ref{DD|Du|^2-1}), an important partial cancellation occurs, and we obtain
\bea \label{DD|Du|^2-2}
F^{p \bar{q}} \na_p \na_{\bar{q}} |\na u|^2_g &=&   - 2 \Re \{ F^{p \bar{q}}  g^{j \bar{i}} u_{\bar{i}} u_p u_{\bar{q} j} \}  + |\na u|_g^2 F^{p \bar{q}} u_{\bar{q} p} -  2 k \Re \{ g^{j \bar{i}} \na_{\bar{i}} u \na_j |\na u|^2_g \} \nonumber\\
&& - 2 \Re \{ g^{j \bar{i}} \E_j u_{\bar{i}} \} + F^{p \bar{q}} g^{j \bar{i}} u_j \hat{R}_{\bar{q} p}{}_{\bar{i}}{}^{\bar{\lambda}} u_{\bar{\lambda}} + |\na \bar{\na} u|^2_{Fg} + |\na \na u|^2_{Fg}.
\eea
We note the identity
\be \label{F-u-id}
F^{p \bar{q}} u_{\bar{q} p} = \Delta_{g} u + \alpha'(k-\gamma) e^{-(1+\gamma)u} a^{p \bar{q}} u_{\bar{q} p} + (k+1) \alpha' \sigma_{k+1} (i \ddb u).
\ee
Substituting the equation (\ref{FY-scalar-g}) into the identity (\ref{F-u-id}), we obtain
\be \label{F-u-id2}
F^{p \bar{q}} u_{\bar{q} p} = - k |\na u|^2_g + \tilde{\E}, 
\ee
where 
\bea \label{tilde-E-def}
\tilde{\E} &=& k \alpha' \sigma_{k+1} (i \ddb u) - \alpha' (k-\gamma)^2 e^{-(1+\gamma)u} a^{p \bar{q}} u_p u_{\bar{q}} \nonumber\\
&& - 2 \alpha' (k-\gamma) e^{-(1+\gamma)u} \Re \{b^i u_i \} - \alpha' e^{-(1+\gamma)u} c + e^{-(k+1)u} \mu.
\eea
will turn out to be another perturbative term. Substituting (\ref{DD|Du|^2-2}) and (\ref{F-u-id2}) into (\ref{grad-main-ineq-1})
\bea \label{grad-main-ineq-2}
F^{p \bar{q}} \na_p \na_{\bar{q}} G &=& {1 \over |\na u|^2_g} |\na \bar{\na} u|^2_{Fg} + {1 \over |\na u|^2_g} |\na \na u|^2_{Fg} -{2 \over |\na u|^2_g} \Re \{ F^{p \bar{q}}  g^{j \bar{i}} u_{\bar{i}} u_p u_{\bar{q} j} \}\nonumber\\
&&- {1 \over |\na u|^4_g} F^{p \bar{q}} \na_p |\na u|^2_g \na_{\bar{q}} |\na u|^2_g  -  2k{1 \over |\na u|^2_g} \Re \{ g^{j \bar{i}} u_{\bar{i}} \na_j |\na u|^2_g \} \nonumber\\
&& - (2+ \sigma) k |\na u|^2_g +{1 \over |\na u|^2_g} F^{p \bar{q}} g^{j \bar{i}} u_j \hat{R}_{\bar{q} p}{}_{\bar{i}}{}^{\bar{\lambda}} u_{\bar{\lambda}}  \nonumber\\
&& -{2 \over |\na u|^2_g} \Re \{ g^{j \bar{i}} \E_j u_{\bar{i}} \} + (2+\sigma) \tilde{\E} .
\eea
Using the critical equation (\ref{critical-eqn}),
\bea \label{apply-crit-eqn}
& \ & - {1 \over |\na u|^4_g} F^{p \bar{q}} \na_p |\na u|^2_g \na_{\bar{q}} |\na u|^2_g  -  2k {1 \over |\na u|^2_g} \Re \{ g^{j \bar{i}} u_{\bar{i}} \na_j |\na u|^2_g \} \nonumber\\
&=& -(1+\sigma)^2 |\na u|^2_F + 2(1+\sigma)k|\na u|^2_g.
\eea
Here we introduced the notation $|\na f|^2_F = F^{p \bar{q}} f_p f_{\bar{q}}$ for a real-valued function $f$. The critical equation (\ref{critical-eqn}) can also be written as 
\be
{ g^{j \bar{i}} \na_p u_j u_{\bar{i}} \over |\na u|_g^2} = -{g^{j \bar{i}} u_j u_{\bar{i} p} \over |\na u|_g^2} - (1+\sigma) u_p.
\ee
We now combine this identity with the Cauchy-Schwarz inequality, which will lead to a partial cancellation of terms. This idea is also used to derive a $C^1$ estimate for the complex Monge-Amp\`ere equation, \cite{Bl,Guan,PS,PSS,Zhang}.
\bea \label{blocki-ineq}
 {1 \over |\na u|^2_g} |\na \na u|^2_{Fg} &\geq& \bigg| {  g^{j \bar{i}} \na u_j u_{\bar{i}} \over |\na u|_g^2}  \bigg|_F^2 \\
&=& {1 \over |\na u|_g^4} |g^{j \bar{i}} u_j \na u_{\bar{i}}|^2_F + (1+\sigma)^2 |\na u|^2_F+ {2(1+\sigma) \over |\na u|^2_g} \Re \{ F^{p \bar{q}} g^{j \bar{i}} u_j u_{\bar{i} p}  u_{\bar{q}} \}.\nonumber
\eea
Let $\varepsilon>0$. Combining (\ref{apply-crit-eqn}) and (\ref{blocki-ineq}) and dropping a nonnegative term,
\bea
& \ & - {1 \over |\na u|^4_g} F^{p \bar{q}} \na_p |\na u|^2_g \na_{\bar{q}} |\na u|^2_g  -  {2k \over |\na u|^2_g} \Re \{ g^{j \bar{i}} u_{\bar{i}} \na_j |\na u|^2_g \}+  (1-\varepsilon) {1 \over |\na u|^2_g} |\na \na u|^2_{Fg} \nonumber\\
&\geq& - (1+\sigma)^2 \varepsilon |\na u|^2_F + 2(1+\sigma)k |\na u|^2_g + {2(1+\sigma)(1-\varepsilon) \over |\na u|^2_g} \Re \{ F^{p \bar{q}} g^{j \bar{i}} u_j u_{\bar{i} p}  u_{\bar{q}} \}.
\eea
Substituting this inequality into (\ref{grad-main-ineq-2}), partial cancellation occurs and we are left with
\bea \label{grad-main-ineq-3}
F^{p \bar{q}} \na_p \na_{\bar{q}} G &\geq& {1 \over |\na u|^2_g} |\na \bar{\na} u|^2_{Fg} + {\varepsilon \over |\na u|^2_g} |\na \na u|^2_{Fg} \nonumber\\
&& +\{ 2\sigma - 2 \varepsilon(1+\sigma) \} {1 \over |\na u|^2_g} \Re \{ F^{p \bar{q}}  g^{j \bar{i}} u_{\bar{i}} u_p u_{\bar{q} j} \}  \nonumber\\
&&+ \sigma k |\na u|^2_g - (1+\sigma)^2 \varepsilon |\na u|^2_F \nonumber\\
&& +{1 \over |\na u|^2_g} F^{p \bar{q}} g^{j \bar{i}} u_j \hat{R}_{\bar{q} p}{}_{\bar{i}}{}^{\bar{\lambda}} u_{\bar{\lambda}}  -{2 \over |\na u|^2_g} \Re \{ g^{j \bar{i}} \E_j u_{\bar{i}} \} + (2+\sigma) \tilde{\E} .
\eea
Since $u \in \Upsilon_k$, we now use (\ref{F-2-g}) in Lemma \ref{F-est-lemma} to pass the norms with respect to $F^{p\bar q}$ to $g^{p\bar q}$ up to an error of order $2^{-6}$. We choose
\be
\varepsilon = (1+\sigma)^{-2} (1+2^{-6})^{-1} {\sigma \over 2}.
\ee
Then
\be
(1+\sigma)^2 \varepsilon |\na u|^2_F \leq {\sigma \over 2} |\na u|^2_g, 
\ee
and
\be
{\varepsilon \over |\na u|^2_g} |\na \na u|^2_{Fg} \geq {\sigma \over 2 (1+\sigma)^2} {1-2^{-6} \over 1+2^{-6}} {1 \over |\na u|^2_g} |\na \na u|^2_{g}.
\ee
Since $\sigma = 2^{-7}$, we have the inequality of numbers ${1 \over 2} {1 - 2^{-6} \over (1+\sigma)^2 (1+2^{-6})} \geq {1 \over 4}$. Thus
\be
{\varepsilon \over |\na u|^2_g} |\na \na u|^2_{Fg} \geq {\sigma \over 4} {1 \over |\na u|^2_g} |\na \na u|^2_{g}.
\ee
We also note the inequalities
\be
{1 \over |\na u|^2_g} |\na \bar{\na} u|^2_{Fg} \geq (1-2^{-6}){1 \over |\na u|^2_g} |\na \bar{\na} u|^2_{g},
\ee
and
\bea
& \ & \{ 2 \sigma - 2 \varepsilon (1+\sigma) \} {1 \over |\na u|^2_g} \Re \{ F^{p \bar{q}}  g^{j \bar{i}} u_{\bar{i}} u_p u_{\bar{q} j} \}  \nonumber\\
&\geq& - \{ 2 - (1+\sigma)^{-1} (1+2^{-6})^{-1} \} \sigma (1+2^{-6}) |\na \bar{\na} u|_{g} \nonumber\\
&\geq& - 2 \sigma (1+2^{-6}) |\na \bar{\na} u|_{g}.
\eea
The main inequality (\ref{grad-main-ineq-3}) becomes
\bea \label{grad-main-ineq-4}
F^{p \bar{q}} \na_p \na_{\bar{q}} G &\geq& (1-2^{-6}) {1 \over |\na u|^2_g} |\na \bar{\na} u|^2_{g} + {\sigma \over 4} {|\na \na u|^2_{g} \over |\na u|^2_g} -2 \sigma(1+2^{-6}) |\na \bar{\na} u|_{g} \nonumber\\
&&+ {\sigma \over 2} |\na u|^2_g +{1 \over |\na u|^2_g} F^{p \bar{q}} g^{j \bar{i}} u_j \hat{R}_{\bar{q} p}{}_{\bar{i}}{}^{\bar{\lambda}} u_{\bar{\lambda}} \nonumber\\
&& -{2 \over |\na u|^2_g} \Re \{ g^{j \bar{i}} \E_j u_{\bar{i}} \} + (2+\sigma) \tilde{\E} .
\eea

\subsection{Estimating the perturbative terms}

\subsubsection{The $\E_j$ terms}
Recall the constant $\Lambda$ is such that $-\Lambda \hat{g}^{j \bar{i}} \leq a^{j \bar{i}} \leq \Lambda \hat{g}^{j \bar{i}}$. We will go through each term in the definition of $\E_j$ $(\ref{E_j})$ and estimate the terms appearing in ${2 \over |\na u|^2_g} \Re \{ g^{j \bar{i}} \E_j u_{\bar{i}} \}$ by groups. In the following, we will use $C$ to denote constants possibly depending on $\alpha'$, $k$, $\gamma$, $a^{p\bar q}$, $b^i$, $c$, $\mu$, and their derivatives.

First, using $2ab \leq a^2 + b^2$, we estimate the terms involving $\nabla\bar\nabla u$
\bea
&&{2|\alpha'(k-\gamma)| \over |\na u|_g^2}e^{-(1+\gamma)u}  | g^{j\bar i} u_{\bar i} (-\gamma  a^{p\bar q} u_{\bar q p} u_j +  \hat\nabla_j a^{p\bar q} u_{\bar q p} +(k-\gamma)  a^{p\bar q} u_p u_{\bar q j} +   \bar b^q u_{\bar q j})|\nonumber\\
&\leq& 2 |\alpha'\Lambda(k-\gamma)(k+2\gamma)| e^{-\gamma} |\nabla \bar\nabla u|_g + C e^{-\gamma u} e^{-u/2} {|\nabla \bar\nabla u|_g\over |\na u|_g} \nonumber\\
&\leq&  2 \bigg\{  |\alpha'\Lambda|^{1/2} (k-\gamma)| \delta^{1/2} |\na u|_g \bigg\} \bigg\{ \delta^{1/2} (k+2\gamma) |\Lambda \alpha'|^{1/2} {|\nabla \bar\nabla u|_g \over |\na u|_g} \bigg\} + C e^{-u/2} {|\nabla \bar\nabla u|_g\over |\na u|_g} \nonumber\\
&\leq&  |\alpha'| \Lambda(k-\gamma)^2 \delta |\na u|^2_g + 4|\Lambda \alpha'|(k+\gamma)^2 \delta  {|\nabla \bar\nabla u|^2_g \over |\na u|^2_g}+ \sigma {|\nabla \bar\nabla u|^2_g \over |\na u|^2_g} + C(\sigma) e^{-u}.
\eea
 
Second, we estimate the terms involving $\nabla\nabla u$ 
\bea
&&{2 |\alpha' (k-\gamma)| \over |\na u|_g^2}e^{-(1+\gamma) u}  | g^{j\bar i} u_{\bar i} \{ (k-\gamma) a^{p\bar q} \nabla_j \nabla_p u \, u_{\bar q}+  b^p \nabla_j \nabla_p u \}|\nonumber\\
&\leq & 2|\alpha'| (k-\gamma)^2 \Lambda e^{-\gamma u} |\nabla\nabla u|_g + 2 \bigg\{ {C \over |\alpha' \Lambda|^{1/2}} e^{-(1+\gamma)u/2} \bigg\} \bigg\{ |\alpha' \Lambda|^{1/2} |k-\gamma| e^{- \gamma u/2} {|\nabla \nabla u|_g \over |\na u|_g} \bigg\} \nonumber\\
&\leq&  |\alpha'| (k-\gamma)^2 \Lambda \delta \bigg\{ {|\nabla \nabla u|^2_g\over |\na u|_g^2} + |\na u|_g^2 \bigg\} + |\alpha' \Lambda| (k-\gamma)^2e^{-\gamma u} {|\nabla \nabla u|^2_g \over |\na u|^2_g} + {C^2 \over  |\alpha' \Lambda|} e^{-(1+\gamma)u} \nonumber\\
&\leq&  2|\alpha'| \Lambda(k-\gamma)^2 \delta {|\nabla \nabla u|^2_g \over |\na u|^2_g} + \delta |\alpha'|(k-\gamma)^2 \Lambda |\na u|_g^2 + C e^{-u}.
\eea
Third, we estimate the terms involving $\nabla u$ quadratically
\bea
&&{2 |\alpha'(k-\gamma)| \over |\na u|_g^2}e^{-(1+\gamma) u}  |g^{j\bar i} u_{\bar i} \{ (k-\gamma) \hat\nabla_j a^{p\bar q} u_p u_{\bar q}- 2(1+\gamma) {\rm Re}\{b^p u_p\} u_j +  u_j b^i u_i \}|
\nonumber\\
&\leq & C e^{-\gamma u} e^{-u/2} |\na u|_g \leq {\sigma\over 16} |\na u|_g^2 + C(\sigma) e^{-(1+2\gamma)u} \leq {\sigma\over 16} |\na u|_g^2 + C e^{-u}.
\eea

Finally, for all the other terms in $\E_j$, we can estimate
\bea
&&{2 |\alpha' (k-\gamma)| \over |\na u|_g^2} e^{-(1+\gamma)u}  |g^{j\bar i} u_{\bar i} \{ -\gamma(k-\gamma) a^{p\bar q} u_p u_{\bar q} u_j + \hat\nabla_j b^p u_p + \partial_j \bar b^q u_{\bar q}  \}|\nonumber\\
&&+{2\over |\na u|_g^2} |g^{j\bar i} u_{\bar i} \{ -(1+\gamma)\alpha' c e^{-(1+\gamma)u} u_j +\alpha' e^{-(1+\gamma)u} \partial_j c + (k+1)e^{-(k+1)u} \mu u_j - e^{-(k+1)u} \partial_j \mu \}| \nonumber\\
&\leq & 
2 |\alpha'| \Lambda (k-\gamma)^2 \gamma  e^{-\gamma u} |\na u|_g^2 + C e^{-(1+\gamma) u} + C e^{-(1+\gamma)u} {e^{-u/2} \over |\na u|_g} + C e^{-(k+1)u} + C e^{-(k+1)u} {e^{-u/2} \over |\na u|_g} \nonumber\\
&\leq & 2 |\alpha'| \Lambda (k-\gamma)^2 \gamma  \delta |\na u|_g^2+ Ce^{-u} + C e^{-u} {e^{-u/2}\over |\na u|_g}. 
\eea

Putting everything together, we obtain the following estimate for the terms coming from $\E_j$.
\bea \label{E_j-estimate}
{2 \over |\na u|^2_g} | g^{j \bar{i}} \E_j u_{\bar{i}} | &\leq& \left\{ 2 |\alpha'|\Lambda (k-\gamma)^2(1+\gamma) \delta  + {\sigma \over 16} \right\} |\na u|^2_g + C e^{-u} + C e^{-u} {e^{-u/2} \over |\na u|_g}  \nonumber
\\
&&+  \{4 |\alpha'| \Lambda (k+\gamma)^2 \delta  + \sigma \} {|\na \bar{\na} u|^2_g \over |\na u|^2_g} +  2 |\alpha'| \Lambda (k-\gamma)^2 \delta  {|\na \na u|^2_g \over |\na u|^2_g} \ \ \ \ \ \ \ \ \ 
\eea

\subsubsection{The $\tilde{\E}$ terms}
Next, estimating $\tilde{\E}$ defined in (\ref{tilde-E-def}) gives
\bea
(2+\sigma) |\tilde{\E}| &\leq& k (2+\sigma)|\alpha'|| \sigma_{k+1}(i \ddb u)| +  (2+\sigma)|\alpha' \Lambda| (k-\gamma)^2 e^{-\gamma u} |\na u|^2_g \nonumber\\
&&+ 2 \|\alpha'(k-\gamma) b^i \|_\infty e^{-\gamma u} e^{-u/2} |\na u|_g + C e^{-(1+\gamma)u} + C e^{-(k+1)u}.
\eea
Using $e^{-u} \leq \delta \leq 1$ and 
\be
2 \|\alpha'(k-\gamma) b \|_\infty e^{-\gamma u} e^{-u/2} |\na u|_g \leq {\sigma \over 16} |\na u|^2_g + C(\sigma) e^{-u} e^{-2 \gamma u},
\ee
we obtain
\bea
(2+\sigma) |\tilde{\E}| \leq k (2+\sigma)|\alpha'| |\sigma_{k+1}(i \ddb u)| +  (2+\sigma)|\alpha' \Lambda| (k-\gamma)^2 \delta  |\na u|^2_g + {\sigma \over 16} |\na u|^2_g + C e^{-u}. \nonumber
\eea
By Lemma \ref{sigma-est-lemma}, we have
\be
k |\alpha'| |\sigma_{k+1}(i \ddb u)| \leq k |\alpha'| {C^{k+1}_n \over n^{1/2} n^{k/2}} |\na \bar{\na} u|_g^k |\na \bar{\na} u|_g \leq \{|\alpha'| C^k_{n-1}  |\na \bar{\na} u|_g^k\} |\na \bar{\na} u|_g.
\ee
Since $u\in \Upsilon_k$, we have $|\alpha'| C^k_{n-1} |\na \bar{\na} u|_g^k \leq 2^{-7}$. Thus
\bea \label{est-tilde-E}
(2+\sigma) |\tilde{\E}| \leq \left\{ (2+\sigma) |\alpha' \Lambda| (k-\gamma)^2 \delta  +  {\sigma \over 16} \right\} |\na u|^2_g +  2^{-7} (2+\sigma)|\na \bar{\na} u|_g +  C e^{-u}.
\eea

\subsection{Completing the estimate}
Combining (\ref{E_j-estimate}) and (\ref{est-tilde-E}),
\bea
{2 \over |\na u|^2_g} | g^{j \bar{i}} \E_j u_{\bar{i}} | + (2+\sigma) |\tilde{\E}| &\leq&  \left\{ 5 |\alpha'|\Lambda (k-\gamma)^2(1+\gamma) \delta  + {\sigma \over 8} \right\} |\na u|^2_g \nonumber\\
&& +  2 |\alpha' \Lambda| (k-\gamma)^2 \delta {|\na \na u|^2_g \over |\na u|^2_g} + \{ 4|\alpha'| \Lambda (k+\gamma)^2 \delta  + \sigma \} {|\na \bar{\na} u|_g^2 \over |\na u|^2_g} \nonumber\\
&&+2^{-7} (2+\sigma) |\na \bar{\na} u|_g  + C e^{-u} + C e^{-u} {e^{-u/2} \over |\na u|_g}.
\eea
Recall $\sigma= 2^{-7}$ and using $(k-\gamma)^2 (1+\gamma) \leq (k+\gamma)^3$, the definition (\ref{delta-def}) of $\delta$ implies 
$$5 |\alpha'|\Lambda (k-\gamma)^2(1+\gamma) \delta  \leq {\sigma \over 8}; \ \ \ \  4|\alpha' \Lambda| (k+\gamma)^2 \delta  \leq 2^{-7}.$$ 
Then, we have
\bea \label{grad-perturb-est}
{2 \over |\na u|^2_g} | g^{j \bar{i}} \E_j u_{\bar{i}} | + (2+\sigma) |\tilde{\E}| &\leq& {\sigma \over 4} |\na u|^2_g +  {\sigma \over 4} {|\na \na u|^2_g \over |\na u|^2_g} + 2^{-6}{ |\na \bar{\na} u|^2_{g} \over |\na u|^2_g} +2^{-7} (2+\sigma) |\na \bar{\na} u|_g  \nonumber\\
&& + C e^{-u} + C e^{-u} {e^{-u/2} \over |\na u|_g}.
\eea
Using (\ref{grad-perturb-est}), the main inequality (\ref{grad-main-ineq-4}) becomes
\bea
F^{p \bar{q}} \na_p \na_{\bar{q}} G &\geq& (1-2^{-5}) {1 \over |\na u|^2_g} |\na \bar{\na} u|^2_{g} - \left\{2\sigma(1+2^{-6}) + 2^{-7} (2+\sigma) \right\} |\na \bar{\na} u|_{g} \nonumber\\
&&+ {\sigma \over 4} |\na u|^2_g +{1 \over |\na u|^2_g} F^{p \bar{q}} g^{j \bar{i}} u_j \hat{R}_{\bar{q} p}{}_{\bar{i}}{}^{\bar{\lambda}} u_{\bar{\lambda}} - C e^{-u} - C e^{-u} {e^{-u/2} \over |\na u|_g} .
\eea
By our choice $\sigma = 2^{-7}$, we have the inequality of numbers
\be
\left\{2\sigma(1+2^{-6}) + 2^{-7} (2+\sigma) \right\}^2 {1 \over 1 - 2^{-5}} \leq {\sigma \over 2}.
\ee
Thus
\bea
& \ & \left\{2\sigma(1+2^{-6})+ 2^{-7} (2+\sigma) \right\} |\na \bar{\na} u|_{g}  \nonumber\\
&\leq&  (1-2^{-5}) {1 \over |\na u|^2_g} |\na \bar{\na} u|^2_{g} + {1 \over 4} \left\{2\sigma(1+2^{-6}) + 2^{-7} (2+\sigma) \right\}^2 {1 \over 1 - 2^{-5}} |\na u|^2_g \nonumber\\
&\leq& (1-2^{-5}) {1 \over |\na u|^2_g} |\na \bar{\na} u|^2_{g} + {\sigma \over 8} |\na u|^2_g.
\eea
We may also estimate 
\be
{1 \over |\na u|^2_g} F^{p \bar{q}} g^{j \bar{i}} u_j \hat{R}_{\bar{q} p}{}_{\bar{i}}{}^{\bar{\lambda}} u_{\bar{\lambda}} \geq - C e^{-u}.
\ee
Putting everything together, at $p$ there holds
\be
0 \geq F^{p \bar{q}} \na_p \na_{\bar{q}} G \geq {\sigma \over 8} |\na u|^2_g - {C e^{-u} e^{-u/2} \over |\na u|_g} - Ce^{-u}.
\ee
From this inequality, we can conclude
\be
|\na u|^2_g(p) \leq C e^{-u(p)}.
\ee
By definition $G(x) \leq G(p)$, and we have
\be
|\na u|^2_g \leq C e^{-u(p)} e^{(1+\sigma)(u(p)-u)} \leq C M^{-1},
\ee
since $e^{u(p)}e^{-u} \leq C$ and $e^{- u} \leq C M^{-1}$ by Theorem \ref{c0-est}. This completes the proof of Theorem \ref{c1-est}.

\section{Second Order Estimate}

\setcounter{equation}{0}

\begin{theorem}
\label{c2-est}
Let $u \in \Upsilon_k$ be a $C^{4}(X)$ function with normalization $\int_X e^u \, \hat{\omega}^n = M$ solving the Fu-Yau equation (\ref{FY-scalar-eq}). Then
\be
| \na \bar{\na} u|^2_g \leq C M^{-1}.
\ee
where $C$ only depends on $(X,\hat{\omega})$, $\alpha'$, $k$, $\gamma$, $\| \rho \|_{C^4(X,\hat{\omega})}$ and $\| \mu \|_{C^2(X)}$. 
\end{theorem}

We begin by noting the following elementary estimate.
\begin{lemma}
Let $\ell \in \{2,3,\dots,n\}$. The following estimate holds:
\be \label{sigma-pqrs-ineq}
|g^{j \bar{i}} \sigma_\ell^{p \bar{q}, r \bar{s}} \na_j u_{\bar{q} p} \na_{\bar{i}} u_{\bar{s} r}| \leq C^{\ell-2}_{n-2} |\na \bar{\na} u|^{\ell-2}_g  |\na \bar{\na} \na u|^2_g.
\ee
\end{lemma}
{\it Proof:}
Since the inequality is invariant, we may work at a point $p \in X$ where $g$ is the identity and $u_{\bar{q} p}$ is diagonal. At $p$, we can use (\ref{secondorder}) and conclude
\be
|g^{j \bar{i}} \sigma_\ell^{p \bar{q}, r \bar{s}} \na_j u_{\bar{q} p} \na_{\bar{i}} u_{\bar{s} r}| \leq \sum_i \sum_{p,q} |\sigma_{\ell -2}(\lambda|pq)| |\na_i u_{\bar{q} p}|^2.
\ee
By Lemma \ref{sigma-est-lemma},
\be
|\sigma_{\ell-2}(\lambda|pq) | \leq {C^{\ell-2}_{n-2} \over (n-2)^{(\ell-2)/2}} |\na \bar{\na} u|^{\ell-2}_g.
\ee
This inequality proves the Lemma. Q.E.D.

\subsection{Differentiating the norm of second derivatives}
\begin{lemma}
Let $u \in \Upsilon_k$ be a $C^4(X)$ function solving (\ref{FY-scalar-eq}) with normalization $\int_X e^u = M$. There exists a constant $C>0$ depending only on $(X,\hat{\omega})$, $\alpha'$, $k$, $\gamma$, $\| \rho \|_{C^4(X,\hat{\omega})}$ and $\| \mu \|_{C^2(X)}$ such that
\bea \label{Fnana-secondorder2}
F^{p \bar{q}} \na_p \na_{\bar{q}} |\na \bar{\na} u|^2_g &\geq& 2 (1-2^{-5}) |\na \bar{\na} \na u|^2_{g} - (1+2k)|\alpha'|^{-1/k} \tau^{1/k} |\na \na u|^2_g \nonumber\\
&& - (1+2k)|\alpha'|^{-1/k} \tau^{1/k} |\na \bar{\na} u|^2_g - |\na \bar{\na} u|_g |\na \bar{\na} \na u|_g \nonumber\\
&& - C M^{-1/2} |\na \bar{\na} \na u|_{g} - CM^{-1} |\na \na u|_g - CM^{-1}.
\eea
\end{lemma}

\medskip

We start by differentiating (\ref{diff-once}) and using the definition of $F^{p \bar{q}}$ to obtain
\bea \label{diff-twice}
0 &=& \alpha' \na_{\bar{i}} \sigma_{k+1}^{p \bar{q}} \na_j \na_p \na_{\bar{q}} u + F^{p \bar{q}} \na_{\bar{i}} \na_j \na_p \na_{\bar{q}} u \nonumber\\
&& + k \na_{\bar{i}}\na_j |\na u|^2_g - \alpha'(k-\gamma)(1+\gamma) e^{-(1+\gamma)u} a^{p \bar{q}} u_{\bar{i}}\na_j \na_p \na_{\bar{q}} u \nonumber\\
&&+ \alpha' (k-\gamma) e^{-(1+\gamma)u} \na_{\bar{i}} a^{p \bar{q}}  \na_j \na_p \na_{\bar{q}} u  + \na_{\bar{i}} \E_j.
\eea
Next, we use (\ref{diff-diff-sigmak}) and (\ref{na-na-na-na-u}) to conclude
\bea \label{diff-twice-2}
F^{p \bar{q}} \na_p \na_{\bar{q}} u_{\bar{i} j} &=& -\alpha' \sigma_{k+1}^{p \bar{q}, r \bar{s}} \na_j \na_p \na_{\bar{q}} u \na_{\bar{i}} \na_r \na_{\bar{s}} u\nonumber\\
&&- F^{p\bar q}\left[ u_p \na_{\bar{i}} \na_j \na_{\bar{q}} u - u_j \na_{\bar{i}} \na_p \na_{\bar{q}} u + u_{\bar{q}} \na_p  \na_{\bar{i}}  \na_j u -  u_{\bar{i}} \na_p \na_{\bar{q}}  \na_j u \right]\nonumber\\
&&- F^{p\bar q}\hat{R}_{\bar{i} j}{}^{\lambda}{}_p u_{\bar{q} \lambda}  +  F^{p\bar q}\hat{R}_{\bar{q} p \bar{i}}{}^{\bar{\lambda}} u_{\bar{\lambda} j} \nonumber\\
&& - k \na_{\bar{i}} \na_j |\na u|^2_g + \alpha' (k-\gamma)(1+\gamma) e^{-(1+\gamma)u} a^{p \bar{q}} u_{\bar{i}}\na_j \na_p \na_{\bar{q}} u \nonumber\\
&&- \alpha'(k-\gamma) e^{-(1+\gamma)u} \na_{\bar{i}} a^{p \bar{q}}  \na_j \na_p \na_{\bar{q}} u  - \na_{\bar{i}} \E_j.
\eea
Direct computation gives
\be \label{direct-nana-u}
F^{p \bar{q}} \na_p \na_{\bar{q}} | \na \bar{\na} u|^2_g = 2 g^{s \bar{i}} g^{j \bar{r}}  F^{p \bar{q}} \na_p \na_{\bar{q}} u_{\bar{i} j} u_{\bar{r} s} + 2 |\na \bar{\na} \na u|^2_{Fgg}.
\ee
Recall (\ref{F-2-g}) that we can pass from $F^{p \bar{q}}$ to the metric $g^{p \bar{q}}$ up to an error of order $2^{-6}$. Substituting (\ref{diff-twice-2}) into (\ref{direct-nana-u}) and estimating terms gives
\bea \label{Fnana-secondorder0}
F^{p \bar{q}} \na_p \na_{\bar{q}} | \na \bar{\na} u|^2_g &\geq&  2 \bigg\{ (1-2^{-6}) |\na \bar{\na} \na u|^2_{g} -  |\alpha' g^{m \bar{i}} g^{j \bar{n}} \sigma_{k+1}^{i \bar{j}, r \bar{s}} \na_j u_{\bar{q} p} \na_{\bar{i}} u_{\bar{s} r} u_{\bar{n} m}| \bigg\} \nonumber\\
&&- C| \na \bar{\na} u|_g |\na \bar{\na} \na u|_g \bigg\{ |D u|_g + e^{-\gamma u} |\na u|_g +  e^{-\gamma u} e^{-{1 \over 2} u}  \bigg\} \nonumber\\
&&- C| \na \bar{\na} u|_g \bigg\{ e^{-u}|\na \bar{\na} u|_g  \bigg\} \nonumber\\
&&-2k  \bigg|  g^{s \bar{i}} g^{j \bar{r}}\na_{\bar{i}}\na_j |\na u|^2_g  u_{\bar{r} s} \bigg| -2 \bigg| g^{s \bar{i}} g^{j \bar{r}}\na_{\bar{i}} \E_j  u_{\bar{r} s} \bigg|.
\eea
The condition $u \in \Upsilon_k$ (\ref{Gamma-cond}) together with $k \leq (n-1)$ gives
\be
C^{k-1}_{n-2} |\alpha'| |\na \bar{\na} u|_g^{k} \leq  |\alpha'| C^k_{n-1} |\na \bar{\na} u|_g^{k} \leq 2^{-7}.
\ee
Therefore by (\ref{sigma-pqrs-ineq})
\be
|\alpha' g^{m \bar{i}} g^{j \bar{n}} \sigma_{k+1}^{p \bar{q}, r \bar{s}} \na_j u_{\bar{q} p} \na_{\bar{k}} u_{\bar{s} r} u_{\bar{n} m}| \leq 2^{-7} |\na \bar{\na} \na u|^2_g.
\ee
In the coming estimates, we will often use the $C^0$ and $C^1$ estimates, and the condition $u \in \Upsilon_k$ (\ref{Gamma-cond}), which we record here for future reference.
\be \label{arsenal}
e^{-u} \leq C M^{-1}, \ \ |\na u|^2_g \leq C M^{-1}, \ \ |\na \bar{\na} u|_g \leq |\alpha'|^{-1/k} \tau^{1/k} ,
\ee
where $\tau = (C^k_{n-1})^{-1} 2^{-7}$. Since $u \in \Upsilon_k$, we have $M=\int_X e^u \hat{\omega}^n \geq 1$, and so we will often only keep the leading power of $M$ since $M \geq 1$. Applying all this to (\ref{Fnana-secondorder0}), we have
\bea \label{Fnana-secondorder}
F^{p \bar{q}} \na_p \na_{\bar{q}}  | \na \bar{\na} u|^2_g  &\geq&  2 (1-2^{-5}) |\na \bar{\na} \na u|^2_{g} \nonumber\\
&& - C M^{-1/2}  | \na \bar{\na} u|_g |\na \bar{\na} \na u|_g - C M^{-1}  | \na \bar{\na} u|_g |\na \bar{\na} u|_g \nonumber\\
&&-2k \bigg|  g^{s \bar{i}} g^{j \bar{r}}\na_{\bar{i}}\na_j |\na u|^2_g  u_{\bar{r} s} \bigg| -2 \bigg| g^{s \bar{i}} g^{j \bar{r}}\na_{\bar{i}} \E_j  u_{\bar{r} s} \bigg|.
\eea
We will now estimate the two last terms. We compute the first of these directly, using (\ref{curv-torsion}) to commute derivatives.
\bea
 2k g^{s \bar{i}} g^{j \bar{r}}\na_{\bar{i}}\na_j |\na u|^2_g  u_{\bar{r} s} &=&  2 k g^{s \bar{i}} g^{j \bar{r}}\bigg\{ g^{p \bar{q}}  u_{\bar{q}} \na_j\na_{\bar{i}} \na_p u + g^{p \bar{q}} u_p \na_{\bar{i}}\na_j \na_{\bar{q}} u \nonumber\\
&& +  g^{p \bar{q}} \na_j \na_p u \na_{\bar{i}}  \na_{\bar{q}} u+ g^{p \bar{q}} u_{\bar{i} p} u_{\bar{q} j} \nonumber\\
&&+  g^{p \bar{q}}  u_{\bar{q}} \hat{R}_{\bar{i} j}{}^\ell{}_p u_\ell  -  g^{p \bar{q}}  u_{\bar{q}} u_{\bar{i} j} u_p \bigg\} u_{\bar{r} s}.
\eea
We estimate
\bea
\bigg|2 k g^{s \bar{i}} g^{j \bar{r}}\na_{\bar{i}}\na_j |\na u|^2_g  u_{\bar{r} s}\bigg| &\leq& k \bigg\{ 4|\na \bar{\na} \na u|_g |\na u|_g + 2|\na \bar{\na} u|_g^2 + 2|\na \na u|_g^2 \nonumber\\
&&+ C e^{-u} |\na u|^2_g + 2|\na u|^2_g |\na \bar{\na} u|_g \bigg\} |\na \bar{\na} u|_g.
\eea
We will use (\ref{arsenal}). Then
\bea \label{Fnana-secondorder-term1}
\bigg|2 k g^{s \bar{i}} g^{j \bar{r}}\na_{\bar{i}}\na_j |\na u|^2_g  u_{\bar{r} s}\bigg| &\leq&  2k |\alpha'|^{-1/k} \tau^{1/k} |\na \bar{\na} u|_g^2+  2k |\alpha'|^{-1/k} \tau^{1/k} |\na \na u|_g^2 \\
&&  + C M^{-1/2} |\na \bar{\na} \na u|_g + CM^{-2} + C M^{-1}. \nonumber
\eea
Next, using the definition (\ref{E_j}) of $\E_j$, we keep track of the order of each term and obtain the estimate
\bea
| g^{s \bar{i}} g^{j \bar{r}} \na_{\bar{i}} \E_j  u_{\bar{r} s}| &\leq& C(a,b,c,\alpha') |\na \bar{\na} u|_g |\na \bar{\na} \na u|_g \bigg\{ e^{- \gamma u} e^{-u/2} + e^{-\gamma u} |\na u|_g  \bigg\} \nonumber\\
&&+C(a,b,c)|\na \bar{\na} u|_g^2 \bigg\{ e^{-\gamma u} |\na u|^2_g + e^{-\gamma u} e^{-u/2} |\na u|_g  + e^{-(1+\gamma)u}  \bigg\} \nonumber\\
&&+C(a,b,c,\alpha') |\na \bar{\na} u|_g  |\na \na u|_g \bigg\{ e^{-\gamma u} |\na u|^2_g + e^{- \gamma u} e^{-u/2} |\na u|_g  + e^{-(1+\gamma)u}  \bigg\} \nonumber\\
&&+C(a,b,c,\alpha') |\na \bar{\na} u|_g \bigg\{e^{-(2+\gamma)u}+ e^{-(1+\gamma) u} e^{-u/2}|\na u|_g+ e^{-(1+\gamma)u} |\na u|^2_g \nonumber\\
&& + e^{-(1+\gamma)u}e^{-u/2} |\na u|^3_g  + e^{-(1+\gamma)u} |\na u|^4_g  \bigg\} \nonumber\\
&&+C(\mu) |\na \bar{\na} u|_g \bigg\{e^{-(k+1)u} |\na \bar{\na} u|_g + e^{-(k+1)u} |\na u|^2_g \nonumber\\
&&+ e^{-(k+1)u} e^{-u/2} |\na u|_g + e^{-(k+2)u} \bigg\} \nonumber\\
&&+  (k-\gamma)^2 g^{s \bar{k}} g^{j \bar{r}} |( \alpha' e^{-(1+\gamma)u} a^{p \bar{q}} \na_{\bar{k}} \na_j \na_p u u_{\bar{q}}) u_{\bar{r} s}| \nonumber\\
&&+  |k-\gamma| g^{s \bar{k}} g^{j \bar{r}} |(\alpha' e^{-(1+\gamma)u} b^i \na_{\bar{k}} \na_j \na_i u)  u_{\bar{r} s}| \nonumber\\
&&+  |k-\gamma| g^{s \bar{k}} g^{j \bar{r}} |(\alpha' e^{-(1+\gamma)u} \gamma a^{p \bar{q}} u_{\bar{q} p} u_{\bar{k} j} u_{\bar{r} s}| \nonumber\\
&&+  (k-\gamma)^2 g^{s \bar{k}} g^{j \bar{r}} |(\alpha' e^{-(1+\gamma)u} a^{p \bar{q}} u_{\bar{k} p} u_{\bar{q}j} )  u_{\bar{r} s}| \nonumber\\
&&+  (k-\gamma)^2 g^{s \bar{k}} g^{j \bar{r}} |(\alpha' e^{-(1+\gamma)u} a^{p \bar{q}} \na_j \na_p u \na_{\bar{k}} \na_{\bar{q}} u)  u_{\bar{r} s}|.
\eea
We will use our estimates (\ref{arsenal}). We also recall the notation $- \Lambda \hat{g}^{p \bar{q}} \leq a^{p \bar{q}} \leq \Lambda \hat{g}^{p \bar{q}}$. We use these estimates and commute covariant derivatives to obtain
\bea
| g^{s \bar{k}} g^{j \bar{r}} \na_{\bar{k}} \E_j u_{\bar{r} s}| &\leq& CM^{-1/2} |\na \bar{\na} \na u|_g + C M^{-1}|\na \na u|_g + C M^{-1} + C M^{-2} \nonumber\\
&&+ CM^{-(k+1)} + CM^{-(k+2)} \nonumber\\
&&+ (k-\gamma)^2   e^{-(1+\gamma) u} g^{s \bar{k}} g^{j \bar{r}} |(\alpha' a^{p \bar{q}} \na_j  \na_{\bar{k}} \na_p u u_{\bar{q}} + \alpha' a^{p \bar{q}} R_{\bar{k} j}{}^\lambda{}_p u_\lambda u_{\bar{q}}) u_{\bar{r} s}| \nonumber\\
&&+ |k-\gamma| e^{-(1+\gamma) u} g^{s \bar{k}} g^{j \bar{r}} |(\alpha' b^i \na_j  \na_{\bar{k}} \na_i u + \alpha' b^i R_{\bar{k} j}{}^\lambda{}_i u_\lambda)  u_{\bar{r} s}| \nonumber\\
&&+2 e^{-\gamma u} |\alpha'| \Lambda (k+\gamma)^2 |\na \bar{\na} u|_g |\na \bar{\na} u|^2_g \nonumber\\
&&+  e^{-\gamma u} |\alpha'| \Lambda(k+\gamma)^2 |\na \bar{\na} u|_g |\na \na u|^2_g.
\eea
Since $u \in \Upsilon_k$, we have $2 |\alpha'| \Lambda (k+\gamma)^2 e^{-\gamma u} \leq 1$.
\bea \label{na-E-j}
| g^{s \bar{k}} g^{j \bar{r}} \na_{\bar{k}} \E_j u_{\bar{r} s}| &\leq& |\alpha'|^{-1/k} \tau^{1/k} |\na \na u|^2_g + |\alpha'|^{-1/k} \tau^{1/k} |\na \bar{\na} u|^2_g \nonumber\\
&&+ C M^{-1/2} |\na \bar{\na} \na u|_g + C M^{-1}|\na \na u|_g + C M^{-1}.
\eea
Substituting (\ref{Fnana-secondorder-term1}) and (\ref{na-E-j}) into (\ref{Fnana-secondorder}) and keeping the leading orders of $M$, we arrive at (\ref{Fnana-secondorder2}).

\subsection{Using a test function}
Let
\be
G = |\na \bar{\na} u|^{2}_g + \Theta |\na u|^2_g.
\ee
where $\Theta \gg 1$ is a large constant depending on $n,k,\alpha'$. To be precise, we let 
\be
\Theta =(1-2^{-6})^{-1}\{ (1+2k)|\alpha'|^{-1/k} \tau^{1/k} + 1\}.
\ee
By (\ref{DD|Du|^2-1}),
\bea
F^{p \bar{q}} \na_p \na_{\bar{q}} |\na u|_g^2 &\geq& |\na \bar{\na} u|^2_{Fg}+ |\na \na u|^2_{Fg} - 2 |\na u|_g |\na \bar{\na} \na u|_g  \nonumber\\
&&- |\na u|^2_g |\na \bar{\na} u|_g - C e^{-u} |\na u|^2_g.
\eea
Applying (\ref{arsenal}) and converting $F^{p \bar{q}}$ to $g^{p \bar{q}}$ yields
\bea  \label{F^pq-|Du|^2-secondorder}
F^{p \bar{q}} \na_p \na_{\bar{q}} |\na u|_g^2 &\geq& (1-2^{-6})|\na \bar{\na} u|^2_{g}+ (1-2^{-6})|\na \na u|^2_{g} - CM^{-1/2} |\na \bar{\na} \na u|_g  \nonumber\\
&&- CM^{-1} |\na \bar{\na} u|_g - C M^{-2}.
\eea
Combining (\ref{Fnana-secondorder2}) and (\ref{F^pq-|Du|^2-secondorder}), we have
\bea \label{Fnana-G}
F^{p \bar{q}} \na_p \na_{\bar{q}} G &\geq&  2 (1-2^{-5}) |\na \bar{\na} \na u|^2_{g} + |\na \bar{\na} u|^2_{g}+  |\na \na u|^2_{g}\nonumber\\
&&-  |\na \bar{\na} u|_g |\na \bar{\na} \na u|_g - C M^{-1/2} |\na \bar{\na} \na u|_{g} \nonumber\\
&& - CM^{-1} |\na \na u|_g - CM^{-1}.
\eea
We will split the linear terms into quadratic terms by applying
\be
CM^{-1/2} |\na \bar{\na} \na u|_g \leq {1 \over 2}  |\na \bar{\na} \na u|_g^2 + {C^2 \over 2 } M^{-1},
\ee
\be
|\na \bar{\na} u|_g |\na \bar{\na} \na u|_g \leq {1 \over 2} |\na \bar{\na} \na u|_g^2  + {1 \over 2} |\na \bar{\na} u|_g^2.
\ee
\be
CM^{-1} |\na \na u|_g \leq  {C^2 \over 4} M^{-2} + |\na \na u|^2_g.
\ee
Applying these estimates, we may discard the remaining quadratic positive terms and (\ref{Fnana-G}) becomes
\be
F^{p \bar{q}} \na_p \na_{\bar{q}} G \geq {1 \over 2} |\na \bar{\na} u|^2_{g} - C M^{-1},
\ee
Let $p \in X$ be a point where $G$ attains its maximum. From the maximum principle, $|\na \bar{\na} u|^2_{g}(p) \leq C M^{-1}$. We conclude from $G \leq G(p)$ that
\be
|\na \bar{\na} u|^2_g \leq C M^{-1}.
\ee
establishing Theorem \ref{c2-est}.

\medskip
We note that many equations involving the derivative of the unknown and/or several Hessians have been studied recently in the literature (see e.g. \cite{CNSI, CNS, CS, GS, GZ, STW, SW, Sun} and references therein). It would be very interesting to determine when estimates with scale hold.

\section{Third Order Estimate}

\setcounter{equation}{0}

\begin{theorem}
\label{c3-est}
Let $u \in {\mathcal{A}}_k$ be a $C^{5}(X)$ function solving equation (\ref{FY-scalar-eq}). Then
\be
| \na \bar{\na} \na u|^2_g \leq C.
\ee
where $C$ only depends on $(X,\hat{\omega})$, $\alpha'$, $k$, $\gamma$, $\| \rho \|_{C^5(X, \hat\omega)}$ and $\|\mu \|_{C^3(X)}$. 
\end{theorem}

To prove this estimate, we will apply the maximum principle to the test function
\be \label{c3-def-G}
G = (|\na \bar{\na} u|^2_g + \eta) |\na \bar{\na} \na u|^2_g + B (|\na u|^2_g + A) |\na \na u|^2_g,
\ee
where $A,B \gg 1$ are large constants to be specified later and $\eta = m \tau^{2/k} |\alpha'|^{-2/k}$. We will specify $m \gg 1$ later and $\tau= (C^k_{n-1})^{-1}2^{-7}$. The condition (\ref{Gamma-cond}) $u \in \Gamma$ implies
\be \label{tau-defn}
|\alpha'|^{1/k} |\na \bar{\na} u|_g \leq \tau^{1/k}.
\ee
Our choice of constants ensures that $\eta$ and $|\na \bar{\na} u|^2_g$ are of the same $\alpha'$ scale. By our previous work, we may estimate by $C$ any term involving $u$, $\na u$, $\na \bar{\na} u$, or the curvature or torsion of $g = e^u \hat{g}$. 
\subsection{Quadratic second order term}

\begin{lemma}
Let $u \in {\mathcal{A}}_k$ be a $C^{4}(X)$ function solving equation (\ref{FY-scalar-eq}). Then for all $A \gg 1$ larger than a fixed constant only depending on $|\na u|_g$ and for all $B>0$, 
\bea \label{quad-second-8}
F^{p \bar{q}} \na_p \na_{\bar{q}} \left\{ ( |\na u|^2_g + A) |\na \na u|^2_g \right\} &\geq& {A \over 2} |\na \na \na u|^2_g + (1-2^{-5})|\na \na u|^4_g\nonumber\\
&&  - {1 \over 2^5 B} |\na \bar{\na} \na u|^4_g - C(A,B).
\eea
where $C(A,B)$ only depends on $A$, $B$, $(X,\hat{\omega})$, $\alpha'$, $k$, $\gamma$, $\| \rho \|_{C^4(X, \hat\omega)}$ and $\|\mu \|_{C^2(X)}$.
\end{lemma}

Differentiating (\ref{diff-once}) gives
\bea \label{quad-second-1}
 F^{p \bar{q}} \na_\ell \na_j \na_p \na_{\bar{q}} u &=& - \alpha' (k-\gamma) \na_\ell ( e^{-(1+\gamma)u} a^{p \bar{q}}) \na_j u_{\bar{q} p}\nonumber\\
&& - \alpha' (\na_\ell \sigma_{k+1}^{p \bar{q}}) \na_j u_{\bar{q} p}  - k \na_\ell \na_j |\na u|^2_g - \na_\ell \E_j,
\eea
Commuting derivatives
\bea \label{quad-second-2}
F^{p \bar{q}} \na_p \na_{\bar{q}} \na_\ell \na_j u&=& F^{p \bar{q}} \na_\ell \na_j \na_p \na_{\bar{q}} u  +  F^{p \bar{q}} \na_p  (\hat{R}_{\bar{q} \ell}{}^\lambda{}_j \na_\lambda u - u_{\bar{q} \ell} u_j) \nonumber\\
&&- F^{p \bar{q}} T^\lambda{}_{p \ell} \na_\lambda \na_j \na_{\bar{q}} u - F^{p \bar{q}} \na_\ell (u_p \na_j \na_{\bar{q}} u - u_j \na_p \na_{\bar{q}} u ).
\eea
We compute directly and commute derivatives to derive
\bea \label{quad-second-3}
F^{p \bar{q}} \na_p \na_{\bar{q}} |\na \na u|^2_g &=& 2 \Re \{ g^{\ell \bar{b}} g^{j \bar{d}} F^{p \bar{q}} \na_p \na_{\bar{q}} \na_\ell \na_j u \na_{\bar{b}} \na_{\bar{d}} u \} \\
&&+  g^{\ell \bar{b}} g^{j \bar{d}} \na_\ell \na_j u  F^{p \bar{q}} R_{\bar{q} p \bar{b}}{}^{\bar{\lambda}} \na_{\bar{\lambda}} \na_{\bar{d}} u+  g^{\ell \bar{b}} g^{j \bar{d}} \na_\ell \na_j u  F^{p \bar{q}} R_{\bar{q} p \bar{d}}{}^{\bar{\lambda}} \na_{\bar{b}} \na_{\bar{\lambda}} u\nonumber\\
&&+  F^{p \bar{q}}  g^{\ell \bar{b}} g^{j \bar{d}} \na_p  \na_\ell \na_j u\na_{\bar{q}} \na_{\bar{b}} \na_{\bar{d}} u +  F^{p \bar{q}}  g^{\ell \bar{b}} g^{j \bar{d}}  \na_{\bar{q}} \na_\ell \na_j u \na_p \na_{\bar{b}} \na_{\bar{d}} u . \nonumber
\eea
Combining (\ref{quad-second-1}), (\ref{quad-second-2}), (\ref{quad-second-3}) and converting $F^{p \bar{q}}$ to $g^{p \bar{q}}$ using Lemma \ref{F-est-lemma}, we estimate
\bea
F^{p \bar{q}} \na_p \na_{\bar{q}} |\na \na u|^2_g &\geq& (1-2^{-6}) |\na \na \na u|^2_g + (1-2^{-6}) |\bar{\na} \na \na u|^2_g \nonumber\\
&&- 2 \alpha' \Re \{ g^{\ell \bar{b}} g^{j \bar{d}} \sigma_{k+1}^{p \bar{q}, r \bar{s}} \na_\ell u_{\bar{s} r} \na_j u_{\bar{q} p} \na_{\bar{b}} \na_{\bar{d}} u \}  - 2 \Re \{ g^{\ell \bar{b}} g^{j \bar{d}} \na_\ell \E_j \na_{\bar{b}} \na_{\bar{d}} u \}   \nonumber\\
&& -C |\na \na u|_g ( |\na \na \na u|_g + |\na \bar{\na} \na u|_g + |\na \na u|_g+ 1) 
\eea
We used the identity (\ref{na-na-na-u}) to estimate $|\na\na\bar\na u|$ by $|\na\bar\na\na u|$ and lower order terms.

Next, using Lemma \ref{sigma-pqrs-ineq} we estimate
\bea
- 2 \Re\{\alpha' g^{\ell \bar{b}} g^{j \bar{d}} \sigma_{k+1}^{p \bar{q}, r \bar{s}} \na_\ell u_{\bar{s} r} \na_j u_{\bar{q} p} \na_{\bar{b}} \na_{\bar{d}} u \} &\geq&  -2 C^{k-1}_{n-2} |\alpha'| |\na \bar{\na} u|_g^{k-1} |\na \na u|_g |\na \bar{\na} \na u|_g^2\nonumber\\
&\geq& -2 C^{k-1}_{n-2} \tau^{1-(1/k)} |\alpha'|^{1/k} |\na \na u|_g |\na \bar{\na} \na u|_g^2 \ \ \  \ 
\eea
and
\be
| g^{\ell \bar{b}} g^{j \bar{d}} \na_\ell \E_j \na_{\bar{b}} \na_{\bar{d}} u  | \leq C|\na \na u|_g \{ 1 + |\na \na u|_g +|\na \bar{\na} \na u|_g +|\na \na \na u|_g \}. 
\ee
Thus
\bea \label{quad-second-4}
F^{p \bar{q}} \na_p \na_{\bar{q}} |\na \na u|^2_g &\geq& (1-2^{-6}) |\na \na \na u|^2_g + (1-2^{-6}) |\bar{\na} \na \na u|^2_g \\
&& -C |\na \na u|_g \{  |\na \bar{\na} \na u|_g^2 + |\na \na \na u|_g + |\na \bar{\na} \na u|_g  + |\na \na u|_g +1\}  \nonumber
\eea
By (\ref{DD|Du|^2-2}),
\be \label{quad-second-5}
F^{p \bar{q}} \na_p \na_{\bar{q}} |\na u|^2_g \geq  (1-2^{-6}) |\na \bar{\na} u|^2_{g} + (1-2^{-6})|\na \na u|^2_{g}  - C |\na \na u|_g - C.
\ee
Direct computation gives
\bea \label{quad-second-6}
F^{p \bar{q}} \na_p \na_{\bar{q}} \left\{ ( |\na u|^2_g + A) |\na \na u|^2_g \right\} &=&( |\na u|^2_g + A)F^{p \bar{q}} \na_p \na_{\bar{q}} |\na \na u|^2_g + |\na \na u|^2_g   F^{p \bar{q}} \na_p \na_{\bar{q}} |\na u|^2_g \nonumber\\
&&+ 2 \Re \{ F^{p \bar{q}} \na_p |\na u|^2_g \na_{\bar{q}} |\na \na u|^2_g \}.
\eea
We estimate
\bea \label{quad-second-7}
2 \left| F^{p \bar{q}} \na_p |\na u|^2_g \na_{\bar{q}} |\na \na u|^2_g \right| 
&\leq& 2(1+2^{-6}) |\na \na u|^2_g |\na u|_g |\bar{\na} \na \na u|_g \nonumber\\
&&+ 2(1+2^{-6}) |\na \na u|^2_g |\na u|_g |\na \na \na u|_g \nonumber\\
&& + C |\bar{\na} \na \na u|_g|\na \na u|_g+C |\na \na \na u|_g|\na \na u|_g.
\eea
Substituting (\ref{quad-second-4}), (\ref{quad-second-5}), (\ref{quad-second-7}) into (\ref{quad-second-6}),
\bea
F^{p \bar{q}} \na_p \na_{\bar{q}} \{ ( |\na u|^2_g + A) |\na \na u|^2_g \} &\geq& A (1-2^{-6}) \left\{ |\na \na \na u|^2_g + |\bar{\na} \na \na u|^2_g \right\}  + (1-2^{-6})|\na \na u|^4_g
\nonumber\\
&&  -3|\na \na u|^2_g |\na u|_g \left\{ |\bar{\na} \na \na u|_g+ |\na \na \na u|_g \right\} \nonumber\\
&& -C(A)  |\na \na u|_g \bigg\{ |\na \bar{\na} \na u|^2_g+ |\na \na \na u|_g + |\na \bar{\na} \na u|_g  \nonumber\\
&& + |\na \na u|^2_g+ |\na \na u|_g + 1 \bigg\} .
\eea
Using $2ab \leq a^2+b^2$,
\be
3|\na \na u|^2_g |\na u|_g |\bar{\na} \na \na u| \leq 2^{-7} |\na \na u|^4_g + 2^5 3^2 |\na u|_g^2 |\bar{\na} \na \na u|_g^2,
\ee
\be
3 |\na \na u|^2_g |\na u|_g |\na \na \na u| \leq  2^{-7} |\na \na u|^4_g + 2^5 3^2 |\na u|_g^2 |\na \na \na u|_g^2,
\ee
\be
C(A) |\na \na \na u|_g |\na \na u|_g \leq |\na \na \na u|_g^2 + {C(A)^2 \over 4} |\na \na u|_g^2
\ee
\be
C(A) |\na \bar{\na} \na u|^2_g |\na \na u|_g \leq {1 \over 2^5 B} |\na \bar{\na} \na u|_g^4 + 2^3 C(A)^2 B |\na \na u|_g^2
\ee
for a constant $B \gg 1$ to be determined later. Then
\bea
F^{p \bar{q}} \na_p \na_{\bar{q}} \left\{ ( |\na u|^2_g + A) |\na \na u|^2_g \right\} &\geq& \left\{ A (1-2^{-6}) -2^6 3^2 |\na u|^2_g -1 \right\} |\na \na \na u|^2_g \nonumber\\
&&+ \left\{ A (1-2^{-6}) -2^6 3^2 |\na u|^2_g -1 \right\} |\bar{\na} \na \na u|^2_g \nonumber\\
&&+ (1-2^{-5})|\na \na u|^4_g  - {1 \over 2^5 B} |\na \bar{\na} \na u|^4_g \\
&&- C(A,B) \bigg\{ |\na \na u|_g+ |\na \na u|^2_g + |\na \na u|^3_g \bigg\}. \nonumber
\eea
The terms $|\na \na u|_g + |\na \na u|^2_g + |\na \na u|^3_g$ can be absorbed into $|\na \na u|_g^4$ by Young's inequality. For $A \gg 1$, obtain (\ref{quad-second-8}).

\subsection{Third order term}

\begin{lemma}
Let $u \in {\mathcal{A}}_k$ be a $C^{5}(X)$ function solving equation (\ref{FY-scalar-eq}). Then
\bea \label{third-order-5}
F^{p \bar{q}} \na_p \na_{\bar{q}} \left\{ (|\na \bar{\na} u|^2_g + \eta) |\na \bar{\na} \na u|^2_g \right\} &\geq&  {1 \over 16} |\na \bar{\na} \na u|^4_g \nonumber\\
&&-C |\na \na \na u|_g \bigg\{ | \na \bar{\na} \na u|_g |\na \na u|_g + | \na \bar{\na} \na u|_g + |\na\na u|_g\bigg\} \nonumber\\
&&- C  \bigg\{| \na \bar{\na} \na u|^2_g |\na \na u|^2_g  + | \na \bar{\na} \na u|^2_g|\na \na u|_g \nonumber\\
&& +  | \na \bar{\na} \na u|_g|\na \na u|^2_g+ | \na \bar{\na} \na u|_g|\na \na u|_g + 1 \bigg\}.
\eea
where $C$ only depends on $(X,\hat{\omega})$, $\alpha'$, $k$, $\gamma$, $\| \rho \|_{C^5(X, \hat\omega)}$ and $\|\mu \|_{C^3(X)}$.
\end{lemma}

To start this computation, we differentiate (\ref{diff-twice-2}).
\bea \label{diff-thrice}
F^{p \bar{q}} \na_i \na_p \na_{\bar{q}} u_{\bar{\ell} j} &=& -\alpha' \na_i (\sigma_{k+1}^{p \bar{q}, r \bar{s}}) \na_j u_{\bar{q} p} \na_{\bar{\ell}} u_{\bar{s} r} -\alpha' \sigma_{k+1}^{p \bar{q}, r \bar{s}} \na_i \na_j u_{\bar{q} p} \na_{\bar{\ell}} u_{\bar{s} r} \nonumber\\
&&-\alpha' \sigma_{k+1}^{p \bar{q}, r \bar{s}}  \na_j u_{\bar{q} p} \na_i \na_{\bar{\ell}} u_{\bar{s} r} + \na_i \left[ -F^{p\bar q} u_p \na_{\bar{\ell}} u_{\bar{q} j} +  F^{p\bar q}u_j \na_{\bar{\ell}} u_{\bar{q}p}  \right] \nonumber\\
&&+ \na_i \left[ -  F^{p\bar q}u_{\bar{q}} \na_p u_{\bar{\ell} j} + F^{p\bar q}u_{\bar{\ell}} \na_p u_{\bar{q} j}  \right] + \na_i [ F^{p\bar q}\hat{R}_{\bar{q} p \bar{\ell}}{}^{\bar{\lambda}} u_{\bar{\lambda} j}- F^{p\bar q}\hat{R}_{\bar{\ell} j}{}^{\lambda}{}_p u_{\bar{q} \lambda} ] \nonumber\\
&& - k \na_i \bigg[ g^{p \bar{q}}  u_{\bar{q}} \na_j u_{\bar{\ell} p}  + g^{p \bar{q}} u_p \na_{\bar{\ell}} u_{\bar{q} j} +  g^{p \bar{q}} \na_j \na_p u \na_{\bar{\ell}}  \na_{\bar{q}} u + g^{p \bar{q}} u_{\bar{\ell} p} u_{\bar{q} j}\nonumber\\
&& +  g^{p \bar{q}}  u_{\bar{q}} \hat{R}_{\bar{\ell} j}{}^\lambda{}_p u_\lambda  -  g^{p \bar{q}}  u_{\bar{q}} u_{\bar{\ell} j} u_p \bigg] + \na_i [ \alpha' (k-\gamma)(1+\gamma) e^{-(1+\gamma)u} a^{p \bar{q}} u_{\bar{\ell}}\na_j u_{\bar{q} p}  ]\nonumber\\
&&- \na_i [\alpha'(k-\gamma) e^{-(1+\gamma)u} \na_{\bar{\ell}} a^{p \bar{q}}  \na_j u_{\bar{q} p} ] - \na_i \na_{\bar{\ell}} \E_j.
\eea
Our conventions (\ref{curv-torsion-defn}) imply the following commutator identities for any tensor $W_{\bar{k} j}$.
\be\label{commu1}
\na_p \na_{\bar{q}} W_{\bar{k} j} = \na_{\bar{q}} \na_p W_{\bar{k} j} + R_{\bar{q} p \bar{k}}{}^{\bar{\lambda}} W_{\bar{\lambda} j} - R_{\bar{q} p}{}^\lambda{}_j W_{\bar{k} \lambda},
\ee
\be\label{commu2}
\na_p \na_{\bar{q}} \na_i W_{\bar{k} j} = \na_i \na_p \na_{\bar{q}} W_{\bar{k} j} + T^\lambda{}_{ip} \na_\lambda W_{\bar{k} j} - \na_p [ R_{\bar{q} i \bar{k}}{}^{\bar{\lambda}} W_{\bar{\lambda} j} - R_{\bar{q} i}{}^\lambda{}_j W_{\bar{k} \lambda}].
\ee
Thus commuting derivatives gives
\bea \label{commute-five}
F^{p \bar{q}} \na_p \na_{\bar{q}} \na_i u_{\bar{k} j} &=&  F^{p \bar{q}} \na_i \na_p \na_{\bar{q}} u_{\bar{k} j} + F^{p \bar{q}} u_i \na_p \na_{\bar{q}} u_{\bar{k} j} - F^{p \bar{q}} u_p \na_i \na_{\bar{q}} u_{\bar{k} j} \nonumber\\
&&+ F^{p \bar{q}} \na_p [ \hat{R}_{\bar{q} i}{}^\lambda{}_j u_{\bar{k} \lambda} - \hat{R}_{\bar{q} i \bar{k}}{}^{\bar{\lambda}} u_{\bar{\lambda} j} ]
\eea
We compute the expression for $F^{p \bar{q}} \na_p \na_{\bar{q}}$ acting on $|\na \bar{\na} \na u|^2_g$, and exchange covariant derivatives to obtain
\bea \label{act-on-nanana}
F^{p \bar{q}} \na_p \na_{\bar{q}} |\na \bar{\na} \na u|^2_g &=& 2 \Re \{g^{i \bar{d}} g^{a \bar{k}} g^{j \bar{b}} F^{p \bar{q}} \na_p \na_{\bar{q}} \na_i u_{\bar{k} j} \na_{\bar{d}} u_{\bar{b} a} \}  \nonumber\\
&& +F^{p \bar{q}} g^{a \bar{d}} g^{e \bar{b}} g^{c \bar{f}} \na_p \na_a u_{\bar{b} c} \na_{\bar{q}} \na_{\bar{d}} u_{\bar{f} e}  +F^{p \bar{q}} g^{a \bar{d}} g^{e \bar{b}} g^{c \bar{f}} \na_a \na_{\bar{q}} u_{\bar{b} c}  \na_{\bar{d}} \na_p u_{\bar{f} e} \nonumber\\
&&  +F^{p \bar{q}} g^{a \bar{d}} g^{e \bar{b}} g^{c \bar{f}} \na_a \na_{\bar{q}} u_{\bar{b} c} R_{\bar{d} p \bar{f}}{}^{\bar{\lambda}} u_{\bar{\lambda} e} -F^{p \bar{q}} g^{a \bar{d}} g^{e \bar{b}} g^{c \bar{f}} \na_a \na_{\bar{q}} u_{\bar{b} c} R_{\bar{d} p}{}^\lambda{}_e u_{\bar{f} \lambda} \nonumber\\
&& -F^{p \bar{q}} g^{a \bar{d}} g^{e \bar{b}} g^{c \bar{f}} R_{\bar{q} a \bar{b}}{}^{\bar{\lambda}} u_{\bar{\lambda} c} \na_p \na_{\bar{d}} u_{\bar{f} e}  +F^{p \bar{q}} g^{a \bar{d}} g^{e \bar{b}} g^{c \bar{f}} R_{\bar{q} a}{}^\lambda{}_c u_{\bar{b} \lambda} \na_p \na_{\bar{d}} u_{\bar{f} e} \nonumber\\
&&+ g^{a \bar{d}} g^{e \bar{b}} g^{c \bar{f}}  \na_a u_{\bar{b} c} F^{p \bar{q}} R_{\bar{q} p \bar{d}}{}^{\bar{\lambda}} \na_{\bar{\lambda}} u_{\bar{f} e} + g^{a \bar{d}} g^{e \bar{b}} g^{c \bar{f}}  \na_a u_{\bar{b} c} F^{p \bar{q}} R_{\bar{q} p \bar{f}}{}^{\bar{\lambda}} \na_{\bar{d}} u_{\bar{\lambda} e} \nonumber\\
&& - g^{a \bar{d}} g^{e \bar{b}} g^{c \bar{f}}  \na_a u_{\bar{b} c} F^{p \bar{q}} R_{\bar{q} p}{}^\lambda{}_e \na_{\bar{d}} u_{\bar{f} \lambda}.
\eea
Substituting (\ref{diff-thrice}) and (\ref{commute-five}) into (\ref{act-on-nanana}), and using Lemma \ref{F-est-lemma}, we have
\bea
F^{p \bar{q}} \na_p \na_{\bar{q}} |\na \bar{\na} \na u|^2_g &\geq& (1-2^{-6}) |\na \na \bar{\na} \na u|^2_g + (1-2^{-6}) |\na \bar{\na} \na \bar{\na} u|^2_g \nonumber\\
&& - 2 \alpha' \Re \{g^{i \bar{d}} g^{a \bar{k}} g^{j \bar{b}} \na_i (\sigma_{k+1}^{p \bar{q}, r \bar{s}}) \na_j u_{\bar{q} p} \na_{\bar{k}} u_{\bar{s} r} \na_{\bar{d}} u_{\bar{b} a} \}\nonumber\\
&& - 2 \alpha' \Re \{g^{i \bar{d}} g^{a \bar{k}} g^{j \bar{b}}  \sigma_{k+1}^{p \bar{q}, r \bar{s}} \na_i \na_j u_{\bar{q} p} \na_{\bar{k}} u_{\bar{s} r}  \na_{\bar{d}} u_{\bar{b} a} \}\nonumber\\
&& - 2 \alpha' \Re \{g^{i \bar{d}} g^{a \bar{k}} g^{j \bar{b}} \sigma_{k+1}^{p \bar{q}, r \bar{s}}  \na_j u_{\bar{q} p} \na_i \na_{\bar{k}} u_{\bar{s} r}  \na_{\bar{d}} u_{\bar{b} a} \}\nonumber\\
&& - C \bigg\{ (|\na \bar{\na} \na \bar{\na} u|_g + |\na \na \bar{\na} \na u|_g) |\na \bar{\na} \na u|_g + |\na \bar{\na} \na \bar{\na} u|_g \nonumber\\
&& + (|\na \na \na u|_g + |\bar{\na} \na \na u|_g +1) |\na \na u|_g |\na \bar{\na} \na u|_g  \nonumber\\
&&+ | \na \bar{\na} \na u|^3_g + | \na \bar{\na} \na u|^2_g + | \na \bar{\na} \na u|_g \bigg\} \nonumber\\
&& -2 \Re \{g^{i \bar{d}} g^{a \bar{k}} g^{j \bar{b}} \na_i \na_{\bar{k}} \E_j \na_{\bar{d}} u_{\bar{b} a} \}
\eea
For the following steps, we will use that $|\alpha'|^{1/k} |\na \bar{\na} u|_g \leq \tau^{1/k}$ for any $u \in {\mathcal{A}}_k$, where $\tau = (C^k_{n-1})^{-1} 2^{-7}$. We also recall that we use the notation $C^\ell_m = {m! \over \ell! (m-\ell)!}$. If $k>1$, we can estimate
\bea
2 | \alpha' g^{i \bar{d}} g^{a \bar{\ell}} g^{j \bar{b}} \na_i (\sigma_{k+1}^{p \bar{q}, r \bar{s}} )\na_j u_{\bar{q} p} \na_{\bar{\ell}} u_{\bar{s} r} \na_{\bar{d}} u_{\bar{b} a} | &\leq& 2 |\alpha'| C^{k-2}_{n-3} |\na \bar{\na} u|^{k-2} |\na \bar{\na} \na u|^4_g \nonumber\\
&\leq&  (2 C^k_{n-1} \tau) |\alpha'|^{2/k} \tau^{-2/k} |\na \bar{\na} \na u|^4_g \nonumber\\
&=& 2^{-6}  |\alpha'|^{2/k} \tau^{-2/k} |\na \bar{\na} \na u|^4_g.
\eea
We used $C^{k-2}_{n-3} \leq C^k_{n-1}$. If $k=1$, the term on the left-hand side vanishes and the inequality still holds. Using the same ideas, we can also estimate
\bea
& \ &  - 2 \alpha' \Re \{g^{i \bar{d}} g^{a \bar{\ell}} g^{j \bar{b}}  \sigma_{k+1}^{p \bar{q}, r \bar{s}} \na_i \na_j u_{\bar{q} p} \na_{\bar{\ell}} u_{\bar{s} r}  \na_{\bar{d}} u_{\bar{b} a} \}  - 2 \alpha' \Re \{g^{i \bar{d}} g^{a \bar{\ell}} g^{j \bar{b}} \sigma_{k+1}^{p \bar{q}, r \bar{s}}  \na_j u_{\bar{q} p} \na_i \na_{\bar{\ell}} u_{\bar{s} r}  \na_{\bar{d}} u_{\bar{b} a} \} \nonumber\\
&\geq& - 2 |\alpha'| C^{k-1}_{n-2} |\na \bar{\na} u|_g^{k-1} |\na \bar{\na} \na u|^2_g \bigg\{ |\na \na \bar{\na} \na u|_g + |\na \bar{\na} \na  \bar{\na} u|_g\bigg\} \nonumber\\
&\geq& - (2 C^{k}_{n-1} \tau) |\alpha'|^{1/k} \tau^{-1/k} |\na \bar{\na} \na u|^2_g \bigg\{ |\na \na \bar{\na} \na u|_g + |\na \bar{\na} \na  \bar{\na} u|_g\bigg\} \nonumber\\
&=&  - 2^{-6} |\alpha'|^{1/k} \tau^{-1/k} |\na \bar{\na} \na u|^2_g \bigg\{ |\na \na \bar{\na} \na u|_g + |\na \bar{\na} \na  \bar{\na} u|_g\bigg\}.
\eea
The perturbative terms can be estimated roughly by using the definition (\ref{E_j}) of $\E_j$ and keeping track of the orders of terms that we do not yet control.
\bea
-2 \Re \{g^{i \bar{d}} g^{a \bar{k}} g^{j \bar{b}} \na_i \na_{\bar{k}} \E_j \na_{\bar{d}} u_{\bar{b} a} \} &\geq& -C |\na \bar{\na} \na u|_g \bigg\{|\na \bar{\na} \na \bar{\na} u|_g  + |\na \bar{\na} \na \na u|_g \nonumber\\
&&+(|\na \bar{\na} \na u|_g+ |\na \na \na u|_g) |\na \na u|_g + |\na \bar{\na} \na u|_g + |\na \na \na u|_g \nonumber\\
&&+|\na \na u|^2_g+ |\na \na u|_g + 1\bigg\}.
\eea
Applying these estimates leads to
\bea \label{third-order-1}
F^{p \bar{q}} \na_p \na_{\bar{q}} |\na \bar{\na} \na u|^2_g &\geq& (1-2^{-6}) \left[ |\na \na \bar{\na} \na u|^2_g + |\na \bar{\na} \na \bar{\na} u|^2_g \right]  - 2^{-6} |\alpha'|^{2/k} \tau^{-2/k} |\na \bar{\na} \na u|^4_g\nonumber\\
&& - 2^{-6} |\alpha'|^{1/k} \tau^{-1/k} |\na \bar{\na} \na u|^2_g \left[  |\na \na \bar{\na} \na u|_g  + |\na \bar{\na} \na  \bar{\na} u|_g \right] \nonumber\\
&&- C {\mathcal{P}}
\eea
where 
\bea \label{defn-P}
{\mathcal{P}} &=& |\na \bar{\na} \na \bar{\na} u|_g |\na \bar{\na} \na u|_g + |\na \na \bar{\na} \na u|_g |\na \bar{\na} \na u|_g +|\na \bar{\na} \na \bar{\na} u|_g \nonumber\\
&&+|\na \na \na u|_g |\na \bar{\na} \na u|_g  |\na \na u|_g + |\na \na \na u|_g |\na \bar{\na} \na u|_g \nonumber\\
&&+ |\na \bar{\na} \na u|_g^2 |\na \na u|_g+ |\na \bar{\na} \na u|_g |\na \na u|^2_g +|\na \bar{\na} \na u|_g |\na \na u|_g\nonumber\\
&& +|\na \na \na u|_g |\na \na u|_g +| \na \bar{\na} \na u|^3_g +| \na \bar{\na} \na u|^2_g  + | \na \bar{\na} \na u|_g.
\eea
We used the fact that the difference between $|\na\bar{\na}\na\na u|_g$ and $|\na\na \bar{\na}\na u|_g$ is a lower order term according to the commutation formula (\ref{commu1}).

Next, we apply (\ref{Fnana-secondorder2}) to obtain
\be \label{third-order-2}
F^{p \bar{q}} \na_p \na_{\bar{q}} |\na \bar{\na} u|^2_g \geq  |\na \bar{\na} \na u|^2_{g} -C |\na \bar{\na} \na u|_g -C |\na \na u|^2_g - C |\na \na u|_g - C.
\ee
We directly compute
\bea \label{third-order-3}
F^{p \bar{q}} \na_p \na_{\bar{q}} \left\{ (|\na \bar{\na} u|^2_g + \eta) |\na \bar{\na} \na u|^2_g \right\} &=& |\na \bar{\na} \na u|^2_g F^{p \bar{q}} \na_p \na_{\bar{q}} |\na \bar{\na} u|^2_g \nonumber\\
&& +  ( |\na \bar{\na} u|^2_g + \eta) F^{p \bar{q}} \na_p \na_{\bar{q}}|\na \bar{\na} \na u|^2_g \nonumber\\
&&+ 2 \Re \{F^{p \bar{q}} \na_p |\na \bar{\na} u|^2_g \na_{\bar{q}} |\na \bar{\na} \na u|^2_g \}.
\eea
We can estimate
\bea \label{third-order-4}
2 \Re \{F^{p \bar{q}} \na_p |\na \bar{\na} u|^2_g \na_{\bar{q}} |\na \bar{\na} \na u|^2_g \} &\geq& -4(1+2^{-6}) |\na \bar{\na} u|_g  |\na \bar{\na} \na u|_g^2 |\na \bar{\na}  \na \bar{\na} u|_g \\
&&-4(1+2^{-6})|\na \bar{\na} u|_g | \na \bar{\na} \na u|_g^2 |\na \na \bar{\na} \na u|_g \nonumber\\
&\geq& -4(1+2^{-6})|\alpha'|^{-1/k} \tau^{1/k} |\na \bar{\na} \na u|_g^2 |\na \bar{\na}  \na \bar{\na} u|_g \nonumber\\
&&-4(1+2^{-6})|\alpha'|^{-1/k} \tau^{1/k}  | \na \bar{\na} \na u|_g^2 |\na \na \bar{\na} \na u|_g . \nonumber
\eea
Combining (\ref{third-order-1}), (\ref{third-order-2}), (\ref{third-order-4}) with (\ref{third-order-3}), setting $\eta = m |\alpha'|^{-2/k} \tau^{2/k}$ and using $|\na \bar{\na} u|_g^2 \leq |\alpha'|^{-2/k} \tau^{2/k}$ leads to
\bea
& \ & F^{p \bar{q}} \na_p \na_{\bar{q}} \left\{ (|\na \bar{\na} u|^2_g + \eta) |\na \bar{\na} \na u|^2_g \right\} \nonumber\\
&\geq& m (1-2^{-6}) |\alpha'|^{-2/k}\tau^{2/k} \bigg\{ |\na \na \bar{\na} \na u|^2_g + |\na \bar{\na} \na \bar{\na} u|^2_g \bigg\} \nonumber\\
&&  -4(1+2^{-6})|\alpha'|^{-1/k} \tau^{1/k} |\na \bar{\na} \na u|_g^2 \bigg\{  |\na \na \bar{\na} \na u|_g + |\na \bar{\na}  \na \bar{\na} u|_g \bigg\} \nonumber\\
&& - 2^{-6} (m+1) |\alpha'|^{-1/k} \tau^{1/k} |\na \bar{\na} \na u|^2_g \bigg\{  |\na \na \bar{\na} \na u|_g +|\na \bar{\na} \na  \bar{\na} u|_g \bigg\} \nonumber\\
&&+ \bigg\{ 1 - 2^{-6} (m+1) \bigg\} |\na \bar{\na} \na u|^4_g-C  |\na \bar{\na} \na u|_g^2|\na \na u|^2_g - C {\mathcal{P}} .
\eea
Using $2ab \leq a^2+b^2$, we estimate
\bea
& \ &  4(1+2^{-6}) |\alpha'|^{-1/k} \tau^{1/k} |\na \bar{\na} \na u|_g^2 \{|\na \bar{\na}  \na \bar{\na} u|_g + |\na \na \bar{\na} \na u|_g \} \nonumber\\
&\leq&  16(1+2^{-6})^2 |\alpha'|^{-2/k} \tau^{2/k} \{ |\na \bar{\na}  \na \bar{\na} u|^2_g +  |\na \na \bar{\na} \na u|^2_g \} + {1 \over 2}  |\na \bar{\na} \na u|_g^4,
\eea
and
\bea
& \ &  2^{-6} (m+1) |\alpha'|^{-1/k} \tau^{1/k} |\na \bar{\na} \na u|^2_g \{ |\na \na \bar{\na} \na u|_g + |\na \bar{\na} \na  \bar{\na} u|_g \}\nonumber\\
&\leq& {1 \over 2} |\alpha'|^{-2/k} \tau^{2/k} \{ |\na \na \bar{\na} \na u|_g^2 + |\na \bar{\na} \na \bar{\na} u|_g^2 \} +  2^{-12} (m+1)^2 |\na \bar{\na} \na u|^4_g .
\eea
The main inequality becomes
\bea
& \ & F^{p \bar{q}} \na_p \na_{\bar{q}} \left\{ (|\na \bar{\na} u|^2_g + \eta) |\na \bar{\na} \na u|^2_g \right\} \nonumber\\
&\geq& \{ m(1-2^{-6}) - 16(1+2^{-6})^2 - {1 \over 2} \} |\alpha'|^{-2/k}\tau^{2/k} \bigg\{ |\na \na \bar{\na} \na u|^2_g + |\na \bar{\na} \na \bar{\na} u|^2_g \bigg\} \nonumber\\
&&+ \left\{ {1 \over 2} - 2^{-6}(m+1) -2^{-12} (m+1)^2  \right\} |\na \bar{\na} \na u|^4_g \nonumber\\
&& -C  |\na \bar{\na} \na u|_g^2|\na \na u|^2_g - C {\mathcal{P}}.
\eea
Next, we estimate terms on the first line in the definition (\ref{defn-P}) of ${\mathcal{P}}$
\bea
& \ & C \{ |\na \bar{\na} \na \bar{\na} u|_g + |\na  \na\bar{\na} \na u|_g \} |\na \bar{\na} \na u|_g \nonumber\\
&\leq& {1 \over 16} |\alpha'|^{-2/k} \tau^{2/k} \{ |\na \bar{\na} \na \bar{\na} u|_g^2 + |\na  \na\bar{\na} \na u|_g^2 \} + 8 C^2 |\alpha'|^{2/k} \tau^{-2/k} |\na \bar{\na} \na u|_g^2
\eea
and
\be
C |\na \bar{\na} \na \bar{\na} u|_g \leq {1 \over 16} |\alpha'|^{-2/k} \tau^{2/k} |\na \bar{\na} \na \bar{\na} u|_g^2 + 4 C^2  |\alpha'|^{2/k} \tau^{-2/k}
\ee
and absorb $| \na \bar{\na} \na u|^3_g + | \na \bar{\na} \na u|^2_g + | \na \bar{\na} \na u|_g $ into $2^{-12} |\na \bar{\na} \na u|^4_g$ plus a large constant. We can now let $m= 18$ and drop the positive fourth order terms. We are left with
\bea
& \ & F^{p \bar{q}} \na_p \na_{\bar{q}} \left\{ (|\na \bar{\na} u|^2_g + \eta) |\na \bar{\na} \na u|^2_g \right\} \nonumber\\
&\geq& \left\{ {1 \over 2} - 2^{-6}(m+1) -2^{-12} (m+1)^2 - 2^{-12} \right\} |\na \bar{\na} \na u|^4_g \nonumber\\
&&-C |\na \na \na u|_g \bigg\{ | \na \bar{\na} \na u|_g |\na \na u|_g + | \na \bar{\na} \na u|_g + |\na\na u|_g\bigg\} \nonumber\\
&&- C  \bigg\{| \na \bar{\na} \na u|^2_g |\na \na u|^2_g  + | \na \bar{\na} \na u|^2_g|\na \na u|_g+ | \na \bar{\na} \na u|_g|\na \na u|^2_g \nonumber\\
&& + | \na \bar{\na} \na u|_g|\na \na u|_g + 1 \bigg\}.
\eea
Since $m=18$,
\be
 {1 \over 2} - 2^{-6}(m+1) -2^{-12} (m+1)^2 - 2^{-12} \geq 2^{-4},
\ee
and we obtain (\ref{third-order-5}).

\subsection{Using the test function}

We have computed $F^{p \bar{q}} \na_p \na_{\bar{q}}$ acting on the two terms of the test function $G$ defined in (\ref{c3-def-G}). Combining (\ref{quad-second-8}) and (\ref{third-order-5})
\bea
F^{p \bar{q}} \na_p \na_{\bar{q}} G &\geq&  {1 \over 32} |\na \bar{\na} \na u|^4_g +{A B \over 2} |\na \na \na u|^2_g + (1-2^{-5})B |\na \na u|^4_g\nonumber\\
&&- C \bigg\{ | \na \na \na u|_g | \na \bar{\na} \na u|_g |\na \na u|_g + |\na \na \na u|_g |\na \bar{\na} \na u|_g + |\na\na\na u|_g |\na\na u|_g\nonumber\\
&&+ | \na \bar{\na} \na u|^2_g |\na \na u|^2_g  + | \na \bar{\na} \na u|^2_g |\na \na u|_g + | \na \bar{\na} \na u|_g |\na \na u|^2_g \nonumber\\
&&+  | \na \bar{\na} \na u|_g |\na \na u|_g \bigg\} - C(A,B). \nonumber
\eea
The negative terms are readily split and absorbed into the positive terms on the first line. For example,
\be
C | \na \na \na u|_g| \na \bar{\na} \na u|_g |\na \na u|_g \leq | \na \na \na u|_g^2 + {C^2 \over 4} | \na \bar{\na} \na u|^2_g |\na \na u|_g^2,
\ee
\be
C | \na \bar{\na} \na u|^2_g |\na \na u|^2_g \leq 2^{-7} |\na \bar{\na} \na u|^4_g + 2^5 C^2 |\na \na u|^4_g 
\ee
\be
C | \na \bar{\na} \na u|^2_g |\na \na u|_g \leq 2^{-7}|\na \bar{\na} \na u|^4_g + 2^5 C^2 |\na \na u|_g^2.
\ee
\be
C | \na \bar{\na} \na u|_g |\na \na u|^2_g \leq 2^{-7} |\na \bar{\na} \na u|^2_g + 2^5 C^2 |\na \na u|_g^4.
\ee
This leads to
\bea
F^{p \bar{q}} \na_p \na_{\bar{q}} G &\geq&  2^{-7} |\na \bar{\na} \na u|^4_g +\{ {A B \over 2} - 1 \} |\na \na \na u|^2_g + \{ {B \over 2} - C \} |\na \na u|^4_g\nonumber\\
&& - C(A,B).
\eea
By choosing $A,B \gg 1$ to be large, we conclude by the maximum principle that at a point $p$ where $G$ attains a maximum, we have 
\be
|\na \bar{\na} \na u|^4_g(p) \leq C, \ \ |\na \na u|^4_g(p) \leq C.
\ee
Therefore $|\na \bar{\na} \na u|_g$ and $|\na\na u|_g$ are both uniformly bounded.

\subsection{Remark on the case $k=1$}

In the case of the standard Fu-Yau equation ($k=1$), to prove Theorem \ref{main-thm} we can instead appeal to a general theorem of concave elliptic PDE and obtain H\"older estimates for the second order derivatives of the solution. To exploit the concave structure, we must rewrite the Fu-Yau equation into the standard form of complex Hessian equation. 

Recall that $\hat{\sigma}_1 (\chi)\, \hat\omega^n = n \chi \wedge \hat\omega^{n-1}, \ \hat{\sigma}_2(\chi) \, \hat{\omega}^n = {n(n-1) \over 2} \chi^2 \wedge \hat{\omega}^{n-2}$. A direct computation with equation (\ref{FY-form-eq}) gives
\bea\label{sigma2fy}
\hat{\sigma}_2(e^u \hat{\omega} + \alpha' e^{-u} \rho + 2\alpha' i\ddb u) &=& 
{n(n-1)\over 2} e^{2u} - 2(n-1) \alpha' e^u |Du|^2_{\hat\omega} - 2(n-1)\alpha' \mu\\
&& + 2(n-1) (\alpha')^2  e^{-u} ( a^{j\bar k} u_j u_{\bar k} - b^i u_i - b^{\bar i} u_{\bar i} )\nonumber\\
&&+ 2(n-1) (\alpha')^2  e^{-u}  c+ (n-1)e^{-u} \hat\sigma_1(\alpha' \rho) + e^{-2u} \hat\sigma_2(\alpha' \rho)\nonumber
\eea
We note that the right hand side of the equation involves the given data $\alpha'$, $\rho$, $\mu$, $u$ and $D u$. 
Since $u\in \Upsilon_1$, the $(1,1)$-form $\o'=e^u \hat{\omega} + \alpha' e^{-u} \rho + 2\alpha' i\ddb u$ is positive definite, and thus both sides of the above equation have a positive lower bound. Moreover, our previous estimates imply that we have uniform a priori estimates on $\| u \|_{C^{1,\beta}(X)}$ for any $0<\beta<1$. The right hand side is therefore bounded in $C^{\beta}(X)$. Since $\hat{\sigma}_2^{1/2}(\chi)$ is a concave uniformly elliptic operator on the space of admissible solutions, we may apply a Evans-Krylov type result of Tosatti-Weinkove-Wang-Yang \cite{TWWY} to conclude $\| u \|_{C^{2,\beta}} \leq C$.

However, for general $k\geq 2$ case, it is impossible to re-write equation (\ref{FY-scalar-eq}) into the standard form of complex Hessian equation and thus there is no obvious concavity we can use.

\bigskip
\noindent
{\it Note:} Just as we were about to post this paper, a preprint,
{\it The Fu-Yau equation in higher dimensions}
by J. Chu, L. Huang, and X.H. Zhu appeared in the net,
arXiv:1801.09351, in which is stated the existence of a solution in the $\Gamma_2$ cone of 
the Fu-Yau equation ($k=1$). Our result is more precise, as our solution is in the
admissible set $\Upsilon_1\subset\Gamma_2$. Moreover, 
our method was used to solve a whole family of Fu-Yau Hessian equations, in which the Fu-Yau equation with $k=1$ is only the simplest example.

\bigskip

Department of Mathematics, Columbia University, New York, NY 10027, USA

\smallskip

phong@math.columbia.edu

\smallskip
Department of Mathematics, Columbia University, New York, NY 10027, USA

\smallskip
 picard@math.columbia.edu

\smallskip
Department of Mathematics, University of California, Irvine, CA 92697, USA

\smallskip
xiangwen@math.uci.edu

\end{document}